\documentclass[twoside,12pt]{article}
\usepackage[latin1]{inputenc}
\usepackage[T1]{fontenc}
\usepackage{amscd}
\usepackage{amsmath}
\usepackage{amsfonts}
\usepackage{amssymb}
\usepackage{amsthm}
\usepackage{comment}

\usepackage[english,francais]{babel}

\author{S. Gou\"ezel}

\newtheorem{thm}{Theorem}[section]
\newtheorem*{thm*}{Theorem}
\newtheorem{prop}[thm]{Proposition}
\newtheorem{lem}[thm]{Lemma}
\newtheorem{cor}[thm]{Corollary}
\newtheorem*{defn}{Definition}
\newtheorem*{notation}{Notation}

\theoremstyle{definition}

\newtheorem*{rmq}{Remark}

\renewcommand{\tilde}{\widetilde}
\renewcommand{\hat}{\widehat}
\newcommand{\norm}[1]{\left\| #1 \right\|}
\newcommand{\dd}{\, {\rm d}}
\newcommand{\N}{\mathbb{N}}
\newcommand{\Z}{\mathbb{Z}}
\newcommand{\R}{\mathbb{R}}
\newcommand{\D}{\mathbb{D}}
\newcommand{\C}{\mathbb{C}}
\newcommand{\tq}{\ |\ }

\newcommand{\boN}{\mathcal{N}}

\renewcommand{\geq}{\geqslant}
\renewcommand{\leq}{\leqslant}
\renewcommand{\phi}{\varphi}
\renewcommand{\epsilon}{\varepsilon}
\newcommand{\boO}{\mathcal{O}}

\newcommand{\boA}{\mathcal{A}}
\newcommand{\boC}{\mathcal{C}}
\DeclareMathOperator{\Hom}{Hom}
\DeclareMathOperator{\real}{Re}
\DeclareMathOperator{\Leb}{Leb}
\newcommand{\psitest}{u}
\newcommand{\moins}{\backslash}

\title{Berry-Esseen theorem and local limit theorem for
non uniformly expanding maps
\footnote{\emph{keywords}: non uniformly expanding maps, Berry-Esseen
theorem, local limit theorem
  \emph{2000 Mathematics Subject Classification:} 37A30, 37A50, 37C30,
  37E05, 47A56, 60F15.}}
\author{S\'ebastien Gou\"ezel
  \footnote{D\'epartement de Math\'ematiques et Applications,
École Normale Sup\'erieure, 45 rue d'Ulm 75005 Paris
(France). e-mail \texttt{Sebastien.Gouezel@ens.fr}}}
\date{October 2004}

\begin{document}
\maketitle

\selectlanguage{english}

\begin{abstract}
In Young towers with sufficiently small tails, the
Birkhoff sums of H\"older continuous functions satisfy a central limit theorem
with speed $O(1/\sqrt{n})$, and a local limit theorem.
This implies the same results for many
non uniformly expanding dynamical systems, namely those
for which a tower with
sufficiently fast returns can be constructed.
\end{abstract}

\selectlanguage{francais}

\begin{abstract}
Dans les tours de Young ayant des queues suffisamment petites, les sommes de
Birkhoff des fonctions h\"olderiennes satisfont le th\'eorème
central limite avec vitesse $O(1/\sqrt{n})$ et le th\'eorème de la
limite locale. Par cons\'equent, de nombreux systèmes dynamiques non
uniform\'ement dilatants satisfont les mêmes conclusions : il suffit
de pouvoir construire une tour avec des retours à la base
suffisamment rapides.
\end{abstract}

\selectlanguage{english}

\section{Results}

\subsection{Introduction}

Let $T:X\to X$ be a probability preserving transformation and $f:X \to
\R$. The functions $f\circ T^k$, for $k\in \N$, are identically
distributed random variables, and it is an important problem in
ergodic theory to see whether they satisfy the same kind of limit
theorems as independent random variables.

Many results are known when $T$ is uniformly expanding or
uniformly hyperbolic (without or with singularities, in the Markov
or non Markov case), and $f$ is H\"older continuous. In this case,
it is indeed often possible to construct a space of functions
containing $f$ on which the transfer operator associated to $T$
has a spectral gap. Therefore, the spectral perturbation method,
introduced by Nagaev in the case of Markov chains, makes it
possible to mimic the probabilistic proofs on independent
variables. In this way, it is possible to get distributional
convergence (to normal laws or stable laws), and more subtle
results such as the speed of convergence (also called the
Berry-Esseen theorem) or the local limit theorem (see for example
\cite{rousseau-egele, guivarch-hardy, broise, aaronson_denker}).
These results, in turn, have important consequences concerning the
asymptotic behavior of the system (\cite{aaronson_denker:2,
varju}).

On the other hand, when the system is not uniformly expanding or uniformly
hyperbolic, it is not possible to use directly the aforementioned
spectral method. Consequently,
other methods have been devised to handle the
distributional convergence of Birkhoff sums. Among many techniques,
the most flexible one is probably the martingale argument of Gordin
(see for example \cite{liverani:CLT, pollicott_sharp, conze_leborgne,
dolgopyat:limit, zweimuller}).
Some results have also been obtained on the speed in the central limit
theorem, by direct estimates
(see \cite{ leborgne_pene, raugi:dim1}). However, there is
currently no result concerning the local limit theorem, which is not
surprising since the proof of
this theorem requires a heavy Fourier machinery, even
in the probabilistic case, and is not easily accessible to elementary
methods.

The aim of this article is to prove the local limit theorem and
the Berry-Esseen theorem for H\"older functions
in the setting of \emph{Young towers}
(\cite{lsyoung:recurrence}), where the decay of correlations is
not exponential and the transfer operator has no spectral gap. The
Young towers are abstract spaces which can be used to model many
non-uniformly expanding maps, for example the
Pomeau-Manneville maps in dimension $1$ studied by Liverani, Saussol
and Vaienti
(\cite{liverani_saussol_vaienti}), the Viana map (for which a
tower is built in \cite{alves_luzzatto_pinheiro}), or the unimodal
maps for which the critical point does not return too quickly
close to itself  (\cite{bruin_luzzatto_vanstrien}). Thus, all
these maps also satisfy the local limit theorem, and the central
limit theorem with speed $O(1/\sqrt{n})$. These results also apply
in non uniformly hyperbolic settings, with the techniques of
\cite{lsyoung:annals}.

The proof is spectral: it uses perturbations of transfer
operators, as in  \cite{guivarch-hardy}, but applied to first
return transfer operators associated to an induced map, as defined by
Sarig in \cite{sarig:decay}. The method
is related to \cite{gouezel:stable}, with a more systematic use of
Banach algebra techniques.

\subsection{Results in Young towers}
\label{section_reduction_tour}

A \emph{Young tower} (\cite{lsyoung:recurrence})
is a probability space $(X,m)$ with a
partition $(B_{i,j})_{i\in I, j<\phi_i}$ of $X$ by positive
measure subsets, where $I$ is finite or countable and
$\phi_i \in \N^*$, together with a nonsingular map $T:X\to X$
satisfying the following properties.
\begin{enumerate}
\item
$\forall i\in I$, $\forall 0\leq j <\phi_i-1$, $T$ is a measure
preserving isomorphism between $B_{i,j}$ and $B_{i,j+1}$.
\item
For every $i\in I$, $T$ is an isomorphism between $B_{i,\phi_i-1}$
and $B:=\bigcup_{k\in I} B_{k,0}$.
\item
Let $\phi$ be the function equal to $\phi_i$ on $B_{i,0}$, whence
$T^\phi$ is a function from $B$ to itself. Let $s(x,y)$ be the
separation time of the points $x$ and $y\in B$ under $T^\phi$,
i.e.\ $s(x,y)=\inf\{ n\tq \exists i\not=j, (T^\phi)^n(x)\in
B_{i,0}, (T^\phi)^n(y)\in B_{j,0}\}$.

As $T^{\phi_i}$ is an isomorphism between $B_{i,0}$ and $B$, it is
possible to consider the inverse $g_m$ of its jacobian with
respect to the measure $m$. We assume that there exist constants
$\beta<1$ and $C>0$ such that $\forall x,y\in B_{i,0}$, $|\log
g_m(x)-\log g_m(y)|\leq C \beta^{s(x,y)}$.
\item
\label{condition_preserve}
The map $T$ preserves the measure $m$.
\item
The partition $\bigvee_{0}^\infty T^{-n}((B_{i,j}))$ separates the
points.
\end{enumerate}

The notion of Young tower has been introduced by Young in
\cite{lsyoung:annals, lsyoung:recurrence} as a model for non
uniformly expanding dynamical systems. The non uniformity is
measured by the \emph{size of tails} $m\{x\in B \tq \phi(x)>n\}$:
if this quantity is very small, then most points enjoy some
expansion before time $n$, when they first return to the basis.
This expansion, in turn, is sufficient to study statistical
properties of the system, including decay of correlations. Young
has proved that, if $m[\phi>n]=O(1/n^\beta)$ for some $\beta>1$,
then the correlations of sufficiently regular functions (see the
definition of $C_\tau(X)$ below) decay like
$O\left(1/n^{\beta-1}\right)$. In particular, if $\beta>2$, these
correlations are summable, and a martingale method can be used to
prove that a central limit theorem holds.

%If we remove the condition \ref{condition_preserve}, Young has shown
%that there exists an absolutely continuous invariant measure.
%Thus, it is not restrictive to assume  \ref{condition_preserve}.

%La condition \ref{condition_pgcd} garantit que $T$ est
%m\'elangeante (sinon, il peut y avoir un problème de p\'eriodicit\'e,
%dont on peut se d\'ebarrasser en regardant une puissance de $T$).
%Encore une fois, on ne perd pas en g\'en\'eralit\'e en supposant cette
%condition.

We extend the separation time $s$ to the whole tower, by setting
$s(x,y)=0$ if $x$ and $y$ are not in the same set $B_{i,j}$, and
$s(x,y)=s(x',y')+1$ otherwise, where $x'$ and $y'$ are the next
iterates of $x$ and $y$ in $B$. For $0<\tau<1$, set
  \begin{equation*}
  C_\tau(X)=\{f: X\to \R \tq \exists C>0 , \forall x,y\in X,
  |f(x)-f(y)|\leq C \tau^{s(x,y)}\}.
  \end{equation*}
This space has a norm $\norm{f}_\tau =\inf\{C
\tq  \forall x,y\in X,
  |f(x)-f(y)|\leq C \tau^{s(x,y)} \} +\norm{f}_\infty$.

The following theorem is well known and can for example be proved
using martingale techniques (see \cite[Theorem 4]{lsyoung:recurrence}).

\begin{thm}
\label{TCL_normal}
Let $\tau<1$. Assume that $m[\phi>n]=O(1/n^{\beta})$ with
$\beta>2$. Let $f \in C_\tau(X)$ have a vanishing integral. Then
there exists $\sigma^2\geq 0$ such that
  \begin{equation*}
  \frac{1}{\sqrt{n}} \sum_{k=0}^{n-1} f\circ T^k
  \to \boN(0,\sigma^2).
  \end{equation*}
Moreover, $\sigma^2=0$ if and only if $f$ is a coboundary, i.e.\
there exists a measurable
function $g$ such that $f=g-g\circ T$ almost everywhere.
\end{thm}

The main results of this article are Theorems
\ref{limite_locale_abstrait} and \ref{TCL_vitesse}. To
formulate the first one, we will need the following definition:

\begin{defn}
A map $f:X \to \R$ is \emph{periodic} if there exist
 $\rho\in \R$, $g:X\to \R$
measurable, $\lambda >0$ and $q:X\to \Z$, such that $f=\rho+g-g\circ
T+\lambda q$ almost everywhere. Otherwise, it is \emph{aperiodic}.
\end{defn}

\begin{thm}[local limit theorem]
\label{limite_locale_abstrait}
Let $\tau<1$. Assume that $m[\phi>n]=O(1/n^{\beta})$ with
$\beta>2$. Let $f \in C_\tau(X)$ have a vanishing integral, and let
$\sigma^2$ be given by Theorem \ref{TCL_normal}.

Assume that $f$ is aperiodic.
This implies in particular
$\sigma^2>0$. Then, for any bounded interval $J\subset \R$, for
any real sequence $k_n$ with $k_n/\sqrt{n} \to \kappa \in \R$, for
any $u\in C_\tau(X)$, for any $v:X\to \R$ measurable,
  \begin{equation*}
  \sqrt{n}\; m\left\{ x \in X \tq S_n f(x) \in J+k_n+u(x)+v(T^n x)\right\}
  \to |J|\frac{e^{-\frac{\kappa^2}{2\sigma^2}}}{\sigma \sqrt{2\pi}}.
  \end{equation*}
\end{thm}

The function on the right is the density of $\boN(0,\sigma^2)$:
this theorem (for $u=v=0$ and $k_n=\kappa \sqrt{n}$) means that
  \begin{equation*}
  m\left\{ \frac{1}{\sqrt{n}} S_n f \in
  \kappa+\frac{J}{\sqrt{n}} \right\} \sim P\left(\boN(0,\sigma^2)
  \in \kappa+\frac{J}{\sqrt{n}} \right).
  \end{equation*}
Hence, it shows that $S_n f/\sqrt{n}$ behaves like $\boN(0,\sigma^2)$
at the local level (contrary to Theorem \ref{TCL_normal} which
deals with the global level).
It is important that $f$ is aperiodic. Otherwise, $f$ could be integer
valued, and the theorem could not hold, e.g.\ for $k_n=0$,
$u=v= 0$ and $J=[1/3,2/3]$.

For $f:X \to \R$, define a function $f_B$ on $B$ by
 \begin{equation}
  \label{definition_fY}
  f_B(x)=\sum_{k=0}^{\phi(x)-1}f(T^k x).
  \end{equation}
In the probabilistic case, the Berry-Esseen theorem, giving the
speed of convergence in the central limit theorem, holds under an $L^3$
moment condition (\cite{feller:2}). In the
dynamical setting, we will need the same kind of hypothesis, but
on the function $f_B$. Note that, since $|f_B|\leq \norm{f}_\infty
\phi$ and $\beta>2$, we always have $f_B \in L^2(B)$.

\begin{thm}[speed in the central limit theorem]
\label{TCL_vitesse} Let $\tau<1$. Assume that
$m[\phi>n]=O(1/n^{\beta})$ with $\beta>2$. Let $f \in C_\tau(X)$
have a vanishing integral, and $\sigma^2$ be given by Theorem
\ref{TCL_normal}.

Assume that $\sigma^2>0$, and that there
exists $0<\delta \leq 1$ such that $\int |f_B|^2 1_{|f_B|>z} \dd
m=O(z^{-\delta})$ when $z\to \infty$. If $\delta=1$, assume also that
$\int f_B^3 1_{|f_B|\leq z}\dd m=O(1)$.
Then there exists $C>0$ such that
$\forall n\in \N^*, \forall a\in \R$,
  \begin{equation*}
  \left| m\left\{ x \tq \frac{1}{\sqrt{n}}S_n f(x) \leq a \right\}
  -P\bigl( \boN(0,\sigma^2) \leq a \bigr) \right|
  \leq \frac{C}{n^{\delta/2}}.
  \end{equation*}
\end{thm}

When $f_B \in L^p$ for some $2<p\leq 3$, then the conditions of the
theorem are satisfied for $\delta=p-2$. In particular, when $f_B \in
L^3$, we obtain a convergence with speed $O(1/\sqrt{n})$, which is the
usual Berry-Esseen theorem.  Note also that, for any $f\in C_\tau(X)$,
the conditions of the theorem are satisfied for $\delta=\beta-2$ if
$2<\beta<3$, and for $\delta=1$ if $\beta>3$. The formulation we have
given is more precise than the usual Berry-Esseen theorem, in view of
the applications, where an $L^p$ condition would not be optimal (see
for example Theorem \ref{thm_LSV}). In fact, the conditions of the
theorem on $f_B$ correspond to necessary and sufficient conditions to
get a central limit theorem with speed $O(n^{-\delta/2})$ in the
probabilistic (independent identically distributed) setting, as shown
in \cite[Theorem 3.4.1]{ibragimov_linnik}.

\begin{rmq}
Using the same methods, it is possible to prove the same results in a
more general setting, namely maps for which a first return map is
\emph{Gibbs-Markov} in the sense of \cite{aaronson:book}. For the sake
of simplicity, we will only consider Young towers.
\end{rmq}

\subsection{Applications}
\label{resultats_metriques}

\subsubsection*{General setting}

Let $(X,d)$ be a locally compact separable metric space, endowed
with a Borel probability measure $\mu$, and $T:X \to X$ a
nonsingular map for which $\mu$ is ergodic.
Assume that there exist a bounded subset $B$ of $X$ with
$\mu(B)>0$, a finite or countable partition (mod $0$) $(B_i)_{i\in
I}$ of $B$, with $\mu(B_i)>0$, and integers $\phi_i
>0$ such that:
\begin{enumerate}
\item $\forall i\in I$, $T^{\phi_i}$ is an isomorphism between $B_i$
and $B$.
\item $\exists \lambda>1$ such that, $\forall i\in I$, $\forall x,y\in
B_i$, $d(T^{\phi_i}x,T^{\phi_i}y) \geq \lambda d(x,y)$.
\item $\exists C>0$ such that, $\forall i\in I$, $\forall x,y \in
B_i$, $\forall k<\phi_i$, $d(T^k x,T^k y) \leq C
d(T^{\phi_i}x,T^{\phi_i}y)$.
\item $\exists \theta>0$ and $D>0$ such that, $\forall i\in I$, the jacobian
$g_\mu$ defined on $B_i$ by $g_\mu(x)=\frac{\dd \mu}{\dd (\mu\circ
T_{|B_i}^{\phi_i} )}$ satisfies: for all $x,y\in B_i$,
$|\log g_\mu(x) -\log g_\mu(y)|\leq D d(T^{\phi_i}x,T^{\phi_i}y)^\theta$.
\end{enumerate}

Denote by $\phi$ the function on $B$ equal to $\phi_i$ on each
$B_i$. If $\mu\{x \tq \phi(x)>n\}$ is summable, we can define a space
$X'=\{ (y,j) \tq y\in B, j<\phi(x)\}$, and a map $T':X' \to X'$ by
$T'(y,j)=(y,j+1)$ if $j<\phi(x)-1$ and $T'(y,j)=(T^{\phi(y)}(y),0)$
otherwise.
Define also $\pi:X' \to X$ by $\pi(y,j)=T^j(y)$. Then $\pi \circ
T'=T\circ \pi$.

Set
$\mu'=\sum_{n=0}^\infty {T'_*}^n( \mu| B\cap \{\phi>n\} )$: it is a
measure of finite mass
on $X'$, not necessarily $T'$-invariant. Young has proved in
\cite[Theorem 1]{lsyoung:recurrence} that there exists a unique
invariant probability measure $m'$ on $X'$ which is absolutely
continuous with respect to $\mu'$. It is ergodic, and $(X',T',\mu')$
is a Young tower in the sense of Section \ref{section_reduction_tour}.
The measure $m=\pi_*(m')$ is
$T$-invariant, absolutely continuous and ergodic.

If $f:X \to \R$ is H\"older continuous, then $f':=f\circ \pi :X' \to
\R$ belongs to $C_\tau(X')$ for $\tau$ close enough to $1$. Moreover,
the Birkhoff sums $\sum_{k=0}^{n-1} f\circ T^k$ and $\sum_{k=0}^{n-1}
f'\circ {T'}^k$ have the same distribution with respect respectively
to $m$ and $m'$. Hence, Theorems \ref{TCL_normal},
\ref{limite_locale_abstrait} and \ref{TCL_vitesse} on the function
$f'$ in the Young tower $(X',T',m')$ imply the same results on the
function $f$ in $(X,T,m)$.

To apply these theorems, we have to check their assumptions. The
condition $m[\phi>n]=O(1/n^\beta)$ with $\beta>2$ corresponds simply
to the requirement
  \begin{equation*}
  \sum_{\phi_i>n} \mu(B_i) = O(1/n^\beta) \text{ for some }\beta>2.
  \end{equation*}

To apply Theorem \ref{limite_locale_abstrait}, we additionally
have to check that the function $f'$ is aperiodic for $T'$, which
can be complicated when the extension $X'$ is not explicitly
described. On the other hand, the aperiodicity of $f$ may be
easier to check, using for example the information at the periodic
points. In this case, the following abstract theorem ensures that
$f'$ is automatically aperiodic, whence we can apply Theorem
\ref{limite_locale_abstrait}.

\begin{thm}
\label{aperiodicite_proj}
Let $T':X'\to X'$ be a probability preserving map on a probability
space $(X',m')$. Let $(X,m)$ be a standard probability space, $T:X\to
X$ an ergodic
probability preserving map, and $\pi:X' \to X$ a map with
countable fibers, such that $m=\pi_*(m')$ and $T\circ \pi=\pi\circ
T'$. Let $f:X\to \R$. Then
\begin{itemize}
\item The function $f$ is a coboundary for $T$ if and only the
function $f\circ \pi$ is a coboundary for $T'$.
\item The function $f$ is aperiodic for $T$ if and only if the
function $f\circ
\pi$ is aperiodic for $T'$.
\end{itemize}
\end{thm}

\subsubsection*{Examples}

Recently, many maps have been shown to fit in the previous setting.
For example,
\cite[Theorem 3]{alves_luzzatto_pinheiro} shows that the Alves-Viana map,
given by
  \begin{equation*}
  T:\left\{ \begin{array}{ccc}
  S^1 \times \R & \to & S^1 \times \R \\
  (\omega,x) & \mapsto & (16\omega, a-x^2+\epsilon \sin(2\pi \omega))
  \end{array}\right.
  \end{equation*}
satisfies these assumptions (for any $\beta>2$) when $0$ is preperiodic
for the map $x\mapsto a -x^2$, and $\epsilon$ is small
enough. In fact, any map close enough to $T$ in the $C^3$-topology
also satisfies them.

In the one-dimensional case, \cite{bruin_luzzatto_vanstrien} shows
that many unimodal maps of the interval also satisfy these
hypotheses: it is sufficient that the returns of the critical point
close to itself occur at a slow enough rate.

Finally, we will discuss with more details the case of the
Pomeau-Manneville maps, studied among many others by
Liverani, Saussol and Vaienti (\cite{liverani_saussol_vaienti}).
They form an interesting class
of applications, since the influence of the
fixed point $0$ becomes more and more
important when $\alpha$ increases. The explicit formula \eqref{eq:LSV}
is not important, what matters is only the local behavior around the
fixed point. Hence, all the following results can be extended to a
much larger class of examples but, for the sake of simplicity, we will
only consider the following maps.

Let $\alpha \in (0,1/2)$, and consider $T:[0,1]\to [0,1]$ given by
  \begin{equation}
  \label{eq:LSV}
   T(x)=\left\{ \begin{array}{cl}
   x(1+2^\alpha x^\alpha) &\text{if }0\leq x\leq 1/2,
   \\
   2x-1 &\text{if }1/2<x\leq 1.
   \end{array}\right.
   \end{equation}
This map has a parabolic fixed point at $0$, and is expanding
elsewhere. It has a unique absolutely continuous invariant
probability measure $m$, whose density is Lipschitz on any
interval of the form $(\epsilon,1]$ (\cite[Lemma
2.3]{liverani_saussol_vaienti}).

\begin{thm}
\label{thm_LSV}
Let $0<\alpha<1/2$, and let $f:[0,1]\to \R$ be a H\"older function
with vanishing integral, which can not be written as $g-g\circ T$.
Then $f$ satisfies a central limit theorem with variance
$\sigma^2>0$.
\begin{itemize}
\item If $\alpha<1/3$, or $f(0)=0$ and there exists
$\gamma>\alpha-1/3$ such that $|f(x)|\leq K x^\gamma$, then
there exists $C>0$ such that $\forall
n\in \N^*$, $\forall a\in \R$,
  \begin{equation*}
  \left| m\left\{ x \tq \frac{1}{\sqrt{n}}S_n f(x) \leq a \right\}
  -P\bigl( \boN(0,\sigma^2) \leq a \bigr) \right|
  \leq \frac{C}{\sqrt{n}}.
  \end{equation*}
\item If $1/3<\alpha<1/2$, $f(0)=0$ and there exists $\gamma>0$ such
that $|f(x)|\leq K x^\gamma$ and $\delta:= \frac{1}{\alpha-\gamma}-2
\in (0,1)$, then there exists $C>0$ such that $\forall
n\in \N^*$, $\forall a\in \R$,
  \begin{equation*}
  \left| m\left\{ x \tq \frac{1}{\sqrt{n}}S_n f(x) \leq a \right\}
  -P\bigl( \boN(0,\sigma^2) \leq a \bigr) \right|
  \leq \frac{C}{n^{\delta/2}}.
  \end{equation*}
\item If $1/3< \alpha<1/2$ and $f(0)\not=0$, then there exists $C>0$
such that $\forall
n\in \N^*$, $\forall a\in \R$,
  \begin{equation*}
  \left| m\left\{ x \tq \frac{1}{\sqrt{n}}S_n f(x) \leq a \right\}
  -P\bigl( \boN(0,\sigma^2) \leq a \bigr) \right|
  \leq \frac{C}{n^{\frac{1}{2\alpha}-1}}.
  \end{equation*}
\end{itemize}
Moreover, if $f$ is aperiodic, it satisfies the local limit theorem.
\end{thm}
\begin{proof}
Let $x_0=1$, and $x_{n+1}$ be the preimage of $x_n$ in $[0,1/2]$.
Let $y_{n+1}$ be the preimage of $x_n$ in $(1/2,1]$: the intervals
$B_n=(y_{n+1}, y_n]$ form a partition of $B=(1/2,1]$ and, if
$\phi_n=n$, all the hypotheses of Section
\ref{resultats_metriques} are satisfied. Moreover, $m(B_n) \sim
\frac{C}{n^{1/\alpha+1}}$ and $x_n \sim \frac{D}{n^{1/\alpha}}$
for constants $C,D>0$ (\cite{liverani_saussol_vaienti}). In
particular, $m(\phi>n)=O(1/n^\beta)$ for $\beta=1/\alpha>2$.

Let $f$ be H\"older on $[0,1]$. If $f(0)\not =0$, then
$f_B=nf(0)+o(n)$ on $B_n$. Otherwise, let $\gamma>0$ be such that
$|f(x)|\leq K x^\gamma$. Reducing $\gamma$ if necessary, we can assume
that $\gamma<\alpha$. Then it is easy to check that $|f_B|\leq C
n^{1-\gamma/\alpha}$ on $B_n$.

Using these estimates, we can check the integrability assumptions of
Theorem \ref{TCL_vitesse} for $\delta=1$ in the first case,
$\frac{1}{\alpha-\gamma}-2$ in the second case, and
$\frac{1}{\alpha}-2$ in the third case. Hence, Theorem
\ref{TCL_vitesse} implies the desired estimates on the speed in the
central limit theorem.

Finally, the local limit theorem is a direct consequence of Theorem
\ref{limite_locale_abstrait}.
\end{proof}

The aperiodicity assumption is \emph{a priori} not easy to check,
since the periodicity equality $f=g-g\circ T+\rho+\lambda q$ is
assumed to hold only almost everywhere. However, under suitable
regularity assumptions on $f$, it is possible to prove that this
equality holds everywhere (see e.g.\ \cite{aaronson_denker} for locally
constant $f$, \cite{gouezel:coboundary} for H\"older $f$). For example, if
$T$ is given by \eqref{eq:LSV}, then $f= \log |T'| - \int \log |T'|$
is aperiodic.

In Section \ref{section_reductions}, we will prove Theorem
\ref{aperiodicite_proj}, and show that it is sufficient to prove
Theorems \ref{limite_locale_abstrait} and \ref{TCL_vitesse} in mixing
Young towers (i.e., such that the return times $\phi_i$ satisfy
$\gcd(\phi_i)=1$). The rest of paper is devoted to the proof of these
theorems.
In Section
\ref{section_resultat_abstrait}, we prove an abstract spectral
result on perturbations of series of operators. In Section
\ref{section_estimees_tour}, we apply this result to first
return transfer operators, to get the key result Theorem
\ref{estimee_clef_markov}. We then use this estimate in the last
two sections to prove respectively the local limit theorem
\ref{limite_locale_abstrait} and the Berry-Esseen theorem
\ref{TCL_vitesse}.

\section{Preliminary reductions}
\label{section_reductions}

\subsection{Proof of Theorem \ref{aperiodicite_proj}}

\begin{proof}[Proof of the coboundary result]
If $f$ is a coboundary, i.e.\ $f=g-g\circ T$, then $f':=f\circ \pi$
can be written as $f'=g'-g'\circ T'$, where $g'=g\circ \pi$. However,
the converse is not immediate: if $f'=g'-g'\circ T'$, the function
$g'$ is \emph{a priori} not constant on the fibers $\pi^{-1}(x)$,
which prevents us from writing $g'=g\circ \pi$.

We use the following characterization of coboundaries: \emph{ Let
$T$ be an endomorphism of a probability space $(X,m)$. Then a
measurable function $f$ on $X$ can be written as $g-g\circ T$ if
and only if}
  \begin{equation}
  \label{cobord_abstrait}
  \forall \epsilon>0, \exists C>0, \forall n\geq 1, m\{x\in X
  \tq |S_n f(x)| \geq C\} \leq \epsilon.
  \end{equation}
This characterization, due to Schmidt, is proved for example in
\cite{aaronson_weiss}.

If $f'$ is a coboundary, then \eqref{cobord_abstrait} is satisfied by
$f'$ in $X'$, whence it is also satisfied by $f$ in $X$ (since this
condition only involves distributions). Thus, $f$ can be written as
$g-g\circ T$.
\end{proof}

\begin{proof}[Proof of the aperiodicity result]
If $f$ is periodic on $X$, i.e.\ $f=\rho+g-g\circ T +\lambda q$
where $q$ is integer-valued, then $f\circ \pi=\rho+(g\circ \pi) -
(g\circ \pi) \circ T' + \lambda (q\circ \pi)$, i.e.\ $f\circ \pi$
is periodic. On the other hand, if $f\circ \pi=\rho'+g'-g'\circ
T'+\lambda q'$, it is not necessarily possible to write directly
$g'=g\circ \pi$. The proof of the
periodicity of $f$ will use ideas of \cite{aaronson_weiss}. We can
assume for example that $\lambda=2\pi$. Replacing $m'$ by one of
its ergodic components, we can also assume that $m'$ is ergodic.

Since the projection $\pi$ has countable fibers, there exists a
measurable subset $A$ of $X'$ such that $\pi$ is an isomorphism
between $A$ and $X$, and $m'(A)>0$. Define a function $\tilde{g}$
on $X$ by $\tilde{g}(x)=g'(x')$, where $x'$ is the unique preimage
of $x$ in $A$. Replacing $f$ by $f-\tilde{g}+\tilde{g}\circ
T-\rho'$, and $g'$ by $g'- \tilde{g}\circ \pi$, we can assume
without loss of generality that $g'=0$ on $A$ and $\rho'=0$.

For $x\in X$, let $W_n(x)$ be the measure on $S^1$ given by
  \begin{equation*}
  W_n(x)=\frac{1}{n}\sum_{k=1}^n \delta(e^{i S_k f(x)})
  \end{equation*}
where $\delta(y)$ is the Dirac mass at $y$. For $u\in C^0(S^1)$,
it is possible by compactness to find a subsequence $n_k$ such that
  \begin{equation*}
  \int_{S^1} u \dd W_{n_k}(x) \to L(u)(x) \text{ weak }\ast
  \text{ in }L^\infty(X).
  \end{equation*}
It is possible to obtain this convergence for a dense countable set
of functions in $C^0(S^1)$, by a diagonal argument. By passing to a
further subsequence, it is also possible to guarantee that
$\frac{1}{n}\sum_{k=1}^n \int_{S^1} u \dd W_{n_k}(x) \to L(u)(x)$ on
a set $Y\subset X$ with $m(Y)=1$, by Komlos' Theorem (\cite{komlos}).
By density, we get the same convergence for any $u \in C^0(S^1)$.

For $x \in Y$, the map $u \in C^0(S^1) \mapsto L(u)(x) \in \R$ is a
nonnegative continuous linear functional sending $1$ to $1$,
thus given by a probability measure
$P_x$. Moreover, these measures satisfy $P_{Tx}(S)=P_x(e^{if(x)} S)$
for any Borel subset $S$ of $S^1$, since
$W_n(Tx)(S)=W_n(x)(e^{if(x)}S) \pm \frac{2}{n}$.

For some $\epsilon>0$, we will prove that
  \begin{equation}
  \label{equation_dirac}
  m\{ x \tq P_x(\{1\})\geq \epsilon \} \geq \epsilon.
  \end{equation}
If $x'\in A\cap {T'}^{-k}(A)$, then $e^{iS_k f\circ \pi(x')} =e^{i
(g'(x')-g'\circ {T'}^k(x'))}=1$, i.e.\ $S_k f\circ \pi(x')\in 2\pi \Z$.
Hence,
  \begin{align*}
  \int_X W_n(x)(\{1\}) \dd m(x)&
  =\frac{1}{n} \sum_{k=1}^n \int_X 1(S_k f(x) \in 2\pi \Z) \dd m(x)
  \\&
  \geq \frac{1}{n}\sum_{k=1}^n \int_{X'} 1(A \cap {T'}^{-k}A) \dd m'(x')
  \\&
  =\int_{X'} 1_A\cdot \left( \frac{1}{n}\sum_{k=1}^n 1_A \circ {T'}^k \right)
  \dd m'(x')
  \to m'(A)^2 >0
  \end{align*}
by Birkhoff Theorem. Thus, for large enough $n$, $\int_X
W_n(x)(\{1\}) \geq 3\epsilon>0$, whence  $\int \left(\frac{1}{n}
\sum_{k=1}^n W_{n_k}(x)\right)(\{1\}) \geq 2\epsilon$ for large
enough $n$. Since $\left(\frac{1}{n} \sum_{k=1}^n
W_{n_k}(x)\right)(\{1\})\leq 1$, we get
  \begin{equation*}
  m\left\{ x \tq \left(\frac{1}{n}
  \sum_{k=1}^n W_{n_k}(x)\right)(\{1\})
  \geq \epsilon \right\}\geq \epsilon.
  \end{equation*}
Thus, the set $C=\left\{x \tq \limsup \left(\frac{1}{n}
  \sum_{k=1}^n W_{n_k}(x)\right)(\{1\}) \geq
\epsilon \right\}$ satisfies $m(C)\geq \epsilon$. Finally,
$P_x(\{1\}) \geq  \epsilon $ on $C$, and this proves
\eqref{equation_dirac}.

Define a measure $\mu$ on $X\times S^1$, by $\mu(U\times V)=\int_U
P_x(V)\dd m(x)$. Then $\mu$ is invariant under the action of
$T_f:(x,y)\mapsto (T(x),e^{-if(x)}y)$, and \cite{aaronson_weiss} proves
that, for almost every ergodic component $P$ of $\mu$, there
exists a compact subgroup $H$ of $S^1$ and a map $\omega:X \to
S^1$ such that, denoting by $m_H$ the Haar measure of $H$,
  \begin{equation*}
  \forall U\times V \subset X \times S^1,\
  P(U\times V)=\int_U m_H(\omega(x)V)\dd m(x).
  \end{equation*}
Moreover, in this case, it is possible to write
$e^{if(x)}=\frac{\omega(Tx)}{\omega(x)} \psi(x)$, where $\psi$ takes
its values in $H$.

If, for all component $P$ of $\mu$, we had $H=S^1$, then $P=m\otimes
\Leb$ for all $P$, whence $\mu=m\otimes \Leb$. This is a
contradiction, since $\mu(C \times \{1\})>0$, while $(m\otimes
\Leb)(C \times \{1\}) =0$. Thus, for some choice of $P$, $H=\Z/k\Z$,
whence $\psi(x)^k=1$, and
$e^{ikf(x)}=\frac{\omega(Tx)^k}{\omega(x)^k}$. Thus, $f$ is periodic
on $X$.
\end{proof}

\begin{rmq}
The proof only shows that the period of $f$ on $X$ divides the
period of $f\circ \pi $ on $X'$, not that they are equal. In fact,
it is not hard to construct examples of Young towers where the two
periods are different. Moreover, the result is not true without
the assumption that the fibers are countable.
\end{rmq}

\subsection{Reduction to the mixing case}

In the proofs in Young towers, it is often useful to assume that
the tower is mixing, i.e.\ $\gcd(\phi_i)=1$. This restriction may
seem technical, but it is important (for example, without it,
there is no decay of correlations any more). For limit theorems,
however, it is irrelevant: \emph{Theorems \ref{TCL_normal},
\ref{limite_locale_abstrait} and \ref{TCL_vitesse} for mixing
towers imply the same results for general towers.}

\begin{proof}
We assume that Theorems \ref{TCL_normal},
\ref{limite_locale_abstrait} and \ref{TCL_vitesse} are true in any
mixing Young tower.
Let $X$ be a non-mixing Young tower, with $N=\gcd(\phi_i)>1$, and let $f\in
C_\tau(X)$ be of vanishing integral.
 For
$k=0,\ldots,N-1$, set $Z_k=\bigcup B_{i,j}$ for $j\equiv k \mod
N$. Then, for every $k$, $(Z_k,T^N)$ is a mixing Young tower, to
which we can apply Theorems  \ref{TCL_normal},
\ref{limite_locale_abstrait} and \ref{TCL_vitesse}.

On $Z_k$, we consider the function $f_k$ given by $\sum_{i=0}^{N-1}
f(T^i x)$, i.e.\ $(S_N f)_{|Z_k}$. Then $\int_{Z_k} f_k=\int_X f=0$.
Theorem \ref{TCL_normal} applied to $f_k$ on $(Z_k,T^N)$ gives a
constant $\sigma_k$ such that, on $Z_k$,
  \begin{equation*}
  \frac{1}{\sqrt{sN}} S_{sN} f \to \boN(0,\sigma_k^2).
  \end{equation*}
Writing an integer $n$ as $sN+r$ with $r<N$, we get that
$\frac{1}{\sqrt{n}}S_n f \to \boN(0,\sigma_k^2)$ on $Z_k$. Finally,
if $x\in Z_0$, $S_n f(x)$ and $S_n f(T^k x)$ differ by at most $2k
\norm{f}_\infty$. Thus, $\frac{1}{\sqrt{n}}S_n f -
\frac{1}{\sqrt{n}}S_n f \circ T^k$ tends to $0$ in probability on
$Z_0$, which shows that $\sigma_k=\sigma_0$. Writing $\sigma$ for
this common number, we get that  $\frac{1}{\sqrt{n}} S_n f \to
\boN(0,\sigma^2)$ on $X$. Moreover, if $\sigma=0$, then $f_0$ is a
coboundary for $T^N$, i.e.\ $f_0=g-g\circ T^N$ where $g:Z_0 \to \R$
is measurable. We extend $g$ to the whole tower: if $(x,j)\in X$,
with $j=sN+r$ and $r<N$, set $g(x,j)=g(x,sN)-\sum_{i=0}^{r-1}
f(x,sN+i)$. It is then easy to check that $f=g-g\circ T$. This proves
Theorem \ref{TCL_normal}.

For the local limit theorem, let us assume that $f$ is aperiodic, and
take $u,v$ and $k_n$ as in the assumptions of Theorem
\ref{limite_locale_abstrait}. We show that the functions $f_k$ are
also aperiodic. Otherwise, for example, $f_0=\rho+ g-g\circ T^N
+\lambda q$, where $g$ and $q$ are defined on $Z_0$. We extend $g$
and $q$ to the whole tower: for $(x,j)\in X$ with $j=sN+r$ and
$0<r<N$, set $g(x,j)=g(x,sN)-\sum_0^{r-1} \left( f(x,sN+i)
-\frac{\rho}{N} \right) +\lambda q(x,sN)$, and $q(x,j)=0$. Then
$f=\frac{\rho}{N}+g-g\circ T+\lambda q$, which is a contradiction.
Thus, all the functions $f_k$ are aperiodic. We can apply to them
Theorem \ref{limite_locale_abstrait} in the mixing tower $(Z_k,T^N)$,
and get that
  \begin{equation*}
  \sqrt{sN}\; m\left\{ x \in Z_k \tq S_{sN} f(x) \in
J+k_{sN}+u(x)+v(T^{sN}x) \right\}
  \to m(Z_k)|J|\frac{e^{-\frac{\kappa^2}{2\sigma^2}}}{\sigma \sqrt{2\pi}}.
  \end{equation*}
Summing over $k$, we get the conclusion of Theorem
\ref{limite_locale_abstrait} for $f$, and for the times of the form
$sN$. For times of the form $sN+r$ with $0<r<N$, we use the same
result for $f\circ T^r$, $u-S_r f, v\circ T^r$
and the sequence $k_{sN+r}$, and get that
  \begin{multline*}
  \sqrt{sN+r}\; m\left\{ x \in X \tq S_{sN}
  f\circ T^r (x) \in J+k_{sN+r}+u(x)-S_r f(x)+v(T^{sN+r}x)\right\}
  \\
  \to |J|\frac{e^{-\frac{\kappa^2}{2\sigma^2}}}{\sigma \sqrt{2\pi}}.
  \end{multline*}
As $S_{sN} f\circ T^r(x) +S_r f(x)=S_{sN+r} f(x)$, this concludes the
proof of Theorem \ref{limite_locale_abstrait}.

Finally, the central limit theorem with speed is deduced from the
same result on each $Z_k$, for times of the form $sN$.  We extend the
result to arbitrary times: writing $n$ as $sN+r$ with $r<N$, we have
$|S_{sN+r}f-S_{sN}f | \leq r \norm{f}_\infty$. This introduces an
error, of order $O(1/\sqrt{n})\leq O(1/n^{\delta/2})$.
\end{proof}

\begin{rmq}
For this proof, it was important to have a strong version of the local
limit theorem, involving functions $u$ and $v$.
\end{rmq}

Theorem \ref{TCL_normal} is proved in
\cite[Theorem 4]{lsyoung:recurrence}.
The rest of the paper is devoted to the proof of Theorems
\ref{limite_locale_abstrait} and
\ref{TCL_vitesse} in mixing Young towers.
\emph{From this point
on, $X$ will be a mixing Young tower, i.e. $\gcd(\phi_i)=1$.\ }

\section{An abstract result}
\label{section_resultat_abstrait}

If $(a_n)_{n\in \N}$ and $(b_n)_{n\in \N}$ are two sequences
indexed by $\N$, we denote by $(a_n)\star (b_n)$ the sequence
$c_n=\sum_{k=0}^n a_k b_{n-k}$. If $(a_n)_{n\in \Z}$ and
$(b_n)_{n\in \Z}$ are two summable sequences indexed by $\Z$, we
also define their convolution $c_n=(a_n)\star (b_n)$ by
$c_n=\sum_{k=-\infty}^{\infty} a_k b_{n-k}$.
Hence, if $\sum a_n z^n$ and $\sum b_n z^n$ are series with summable
coefficients, the coefficient of $z^n$ in $\left(\sum a_k z^k\right)
\left(\sum b_k z^k \right)$ is given by $(a_n)\star (b_n)$ (more
precisely, it is given by the $n^{\text{th}}$ term $\bigl((a_k)\star
(b_k)\bigr)_n$ of this sequence, but we will often abuse notations and
write simply $(a_n)\star (b_n)$).
Finally, we write
 $\D=\{z \in \C \tq |z|<1\}$ and $\overline{\D}=\{z\in \C \tq
|z|\leq 1\}$.

The goal of this section is to prove the following theorem:
\begin{thm}
\label{thm_abstrait_vitesse}
Let $\beta>2$. Let $R_n$, for $n\in \N^*$, be operators on a
Banach space $E$, with $\sum_{k=n+1}^\infty \norm{R_k}=O(1/n^{\beta})$. Set
$R(z)=\sum R_n z^n$, and assume that $1$ is a simple isolated
eigenvalue of $R(1)$, while $I-R(z)$ is invertible for $z \in
\overline{\D}-\{1\}$. Let $P$ be the spectral projection
associated to $R(1)$ and the eigenvalue $1$, and assume that
$PR'(1)P=\mu P$ with $\mu>0$.

Let $R_n(t)$ be operators on $E$ (for $t$ in some interval
$[-\alpha,\alpha]$ with $\alpha>0$) such that $\sum_{k=n+1}^\infty
\norm{R_k(t)-R_k}
\leq C |t|/n^{\beta-1}$. Set $R(z,t)=\sum z^n R_n(t)$. Let
$\lambda(z,t)$ be the eigenvalue close to $1$ of $R(z,t)$, for
$(z,t)$ close to $(1,0)$. We assume that $\lambda(1,t)=1-M(t)$
with $M(t)\sim c t^2$ for some constant $c$ with $\real(c)>0$.

Then, for small enough $t$, $I-R(z,t)$ is invertible for all
$z\in\D$. Let us denote its inverse by $\sum T_{n}(t)z^n$. Then there
exist $\alpha'>0$, $d>0$ and $C>0$ such that, for every $t\in
[-\alpha',\alpha']$, for every $n\in \N^*$,
  \begin{equation}
  \label{EquationTCLabstrait}
  \norm{T_{n}(t)-\frac{1}{\mu}\left(1-\frac{1}{\mu}M(t)\right)^n P }
  \leq \frac{C}{n^{\beta-1}}
  +C|t| \left( \frac{1}{n^{\beta-1}} \right) \star \left(1-dt^2\right)^n.
  \end{equation}
\end{thm}

In the application of Theorem \ref{thm_abstrait_vitesse} to the
proof of Theorems \ref{limite_locale_abstrait} and
\ref{TCL_vitesse}, the operators $R_n$ will describe the returns
to the basis $B$, and will be easily understood, as well as their
perturbations $R_n(t)$. On the other hand, $T_n$ will describe all
the iterates at time $n$, and $T_{n}(t)$ will be closely related to
the characteristic function $E(e^{it S_n f})$. Thus,
\eqref{EquationTCLabstrait} will enable us to describe precisely
$E(e^{it S_n f})$, and this information will be sufficient to get
Theorems \ref{limite_locale_abstrait} and \ref{TCL_vitesse}.

\subsection{Banach algebras and Wiener Lemma}

\label{subsection_wiener} In this paragraph, we define some Banach
algebras which will be useful in the following estimates. We postpone
the proofs of the properties of these algebras to Appendix
\ref{appendice_preuve_wiener}.

A Banach algebra $\boA$ is a complex Banach space with an associative
multiplication $\boA\times \boA \to \boA$ such that $\norm{AB}\leq
\norm{A}\norm{B}$, and a neutral element. The set of invertible
elements is then an open subset of $\boA$, on which the inversion is
continuous.

Let $\boC$ be a Banach algebra. If $\gamma>1$, we write
$\boO_\gamma( \boC)$ for the set of formal series
$\sum_{n=-\infty}^\infty A_n z^n$, where $A_n\in \boC$ and
$\norm{A_n}=O(1/|n|^\gamma)$ when $n \to \pm\infty$, endowed with the
standard product of power series, corresponding to the convolution of
the sequences $(A_n)$ and $(B_n)$. It admits the
naive norm $ \sup_{n\in \Z} (|n|+1)^{\gamma}\norm{A_n}$, for which it is not
a Banach algebra. However, there exists a norm, equivalent to the
previous one, which makes $\boO_\gamma(\boC)$ into a Banach algebra
(Proposition \ref{wiener_algebra}).
 Moreover, this algebra satisfies a Wiener
Lemma: if $A(z)=\sum A_n z^n \in \boO_\gamma(\boC)$ is such that
$A(z)$ is invertible for every $z\in S^1$, then $A$ is
invertible in $\boO_\gamma(\boC)$ (Theorem \ref{wiener_S1}).

We will also use the Banach algebra $\boO_\gamma^+(\boC)$, given
by the set of series $\sum A_n z^n \in \boO_\gamma(\boC)$ such
that $A_n=0$ for $n<0$. It is a closed subalgebra of
$\boO_\gamma(\boC)$, and it satisfies also a Wiener Lemma (Theorem
\ref{wiener_D}).

\begin{notation}
If $f: [-\alpha,\alpha] \times \Z\to \R_+$ for some $\alpha>0$, and
$\boC$ is a Banach algebra, we denote by $O_\boC(f(t,n))$ the set of
series $\sum_{n=-\infty}^\infty c_n(t) z^n$ where
$c_n:[-\alpha,\alpha]\to \boC$ is such that there exists $\alpha'>0$
and $C>0$ such that
  \begin{equation*}
  \forall t\in [-\alpha',\alpha'], \forall n \in \Z,\
  \norm{c_n(t)} \leq C f(t,n).
  \end{equation*}
\end{notation}
We will often omit the subscript in $O_\boC$. As usual, we will often
write $\sum c_n(t)z^n = O(f(t,n))$ instead of the more correct
formulation $\sum c_n(t) z^n \in O(f(t,n))$.
We will also write
$O(g(n))$ for the set of series $\sum A_n z^n$ with
$\norm{A_n}\leq C g(n)$ for some constant $C$. This is a particular
case of the previous notation, where the functions $f(t,n)$ are
independent of $t$. Until the end of Section
\ref{section_resultat_abstrait}, the notation $O$ will always have
this signification.

\begin{rmq}
The notations $\boO$ and $O$ should not be confused: there are
similarities (which is why we have used the same letter), but the
calligraphic notation $\boO$ indicates additionally a Banach algebra. In
this case, we can for example use the continuity of inversion.
\end{rmq}

With these notations, we can reformulate Theorem \ref{wiener_S1}
as follows: if $\sum A_n z^n =O(1/(|n|+1)^\gamma)$ for $\gamma>1$, and
$\sum A_n z^n$ is invertible for every $z\in S^1$, then $\left(
\sum A_n z^n \right)^{-1} = O(1/(|n|+1)^\gamma)$. The fact that
$\boO_\gamma(\boC)$ is a Banach algebra also implies that, for
$\gamma>1$,
  \begin{equation}
  O\left(\frac{1}{(|n|+1)^\gamma}\right) \star
  O\left(\frac{1}{(|n|+1)^\gamma}\right)
  \subset O\left(\frac{1}{(|n|+1)^\gamma}\right),
  \end{equation}
i.e., if two series $\sum A_n z^n$ and $\sum B_n z^n$ (with $A_n,B_n
\in \boC$) satisfy $\sup_{n\in \Z} (|n|+1)^\gamma \norm{A_n} <\infty$
and $\sup_{n\in \Z} (|n|+1)^\gamma \norm{B_n}<\infty$, then the series
$\sum C_n z^n:= \left(\sum A_n z^n\right)\left(\sum B_n z^n \right)$
also satisfies $\sup_{n\in \Z} (|n|+1)^\gamma \norm{C_n}<\infty$.

\subsection{Preliminary technical estimates}

For notational convenience, we will often write $t$ instead of $|t|$
in what follows. Equivalently, the reader may consider that the
proofs are written for $t\geq 0$. We will also write
$\frac{1}{|n|^\gamma}$ instead of $\frac{1}{(|n|+1)^\gamma}$,
discarding the problem at $n=0$.

\begin{lem}
\label{ngamma_star_t}
When $\gamma>1$ and $d>0$,
  \begin{equation*}
  O\left(\frac{1}{|n|^\gamma}\right)\star
  O\left(1_{n\geq 0} t^2 (1-dt^2)^n\right)
  \subset
  O\left(\frac{1}{|n|^\gamma} + 1_{n\geq 0} t^2
  \left(1-\frac{d}{2}t^2\right)^n \right).
  \end{equation*}
\end{lem}
\begin{proof}
For $n<0$, the coefficient in the convolution is less than
$\sum_{k=0}^\infty t^2(1-dt^2)^k \frac{1}{|n|^\gamma} \leq
\frac{C}{|n|^\gamma}$. For $n\geq 0$, it is less than
$\sum_{k=-\infty}^{n/2} \frac{1}{|k|^\gamma}
t^2(1-dt^2)^{n/2}+\sum_{k=n/2}^n \frac{1}{(n/2)^\gamma}
t^2(1-dt^2)^{n-k} \leq C t^2(1-dt^2)^{n/2}+\frac{C}{n^\gamma}$.
Finally, as $\sqrt{1-dt^2}\leq \left(1-\frac{d}{2}t^2\right)$, we get
the conclusion.
\end{proof}

\begin{lem}
\label{DL_inversion} Let $\gamma>1$ and $d>0$. Let
$G_t(z)=O\left(\frac{t}{|n|^\gamma}+1_{n\geq 0}
t^3(1-dt^2)^{n}\right)$, and assume that $F(z)=O(1/|n|^\gamma)$ is
invertible for every $z\in S^1$. Then $\bigl[F(z)+G_t(z)
\bigr]^{-1}=F(z)^{-1}+O\left(\frac{t}{|n|^\gamma}+1_{n\geq 0}t^3
\left(1-\frac{d}{64}t^2\right)^n\right)$.
\end{lem}
\begin{proof}
We first assume that $F(z)=1$. Setting $H_t(z)=\sum_{n\in \Z}
\frac{t}{|n|^\gamma}z^n + \sum_{n\in \N} t^3 (1-dt^2)^nz^n$, the
norm of the coefficients of $[1+G_t(z)]^{-1}$ is less than the
coefficients of $[1-H_t(z)]^{-1}$. Thus, it is sufficient to
consider $\frac{1}{1-tK(z)-\frac{t^3}{1-(1-dt^2)z}}$ where
$K(z)=\sum \frac{z^n}{|n|^\gamma}$. Note that
  \begin{equation}
  \label{yosjadgh}
  \frac{1}{1-tK(z)-\frac{t^3}{1-(1-dt^2)z}}
  =\frac{1}{1-tK(z)-t^3}
  \frac{1-(1-dt^2)z}{1-(1-dt^2)\left[1+\frac{t^3}{1-tK(z)-t^3}\right]z}.
  \end{equation}
For small enough $t$,
$\left|(1-dt^2)\left[1+\frac{t^3}{1-tK(z)-t^3}\right]\right| < 1$,
whence
  \begin{equation}
  \label{aerj;sdf}
  \begin{split}
  \frac{1-(1-dt^2)z}{1-(1-dt^2)\left[1+\frac{t^3}{1-tK(z)-t^3}\right]z}
  \!\!\!\!\!\!\!\!\!\!\!\!\!\!\!\!\!\!\!\!\!\!\!\!
  \!\!\!\!\!\!\!\!\!\!\!&
  \\&
  =\left(1-(1-dt^2)z\right)\sum_{n=0}^{\infty} (1-dt^2)^n
  \left[1+\frac{t^3}{1-tK(z)-t^3}\right]^n z^n
  \\&
  =1+\frac{t z (1-dt^2)}{1-tK(z)-t^3}t^2\sum_{n=0}^\infty (1-dt^2)^n
  \left[1+\frac{t^3}{1-tK(z)-t^3}\right]^n z^n.
  \end{split}
  \end{equation}
We first study the sum. Let $\boA=\boO_\gamma(\C)$ be the Banach
algebra of the series whose coefficients are $O(1/|n|^\gamma)$. As
$K(z)\in \boA$, we have $1-tK(z)-t^3=1+O_\boA(t)$. Since the
inversion is Lipschitz on a Banach algebra, we get
$\frac{t^3}{1-tK(z)-t^3}=t^3+O_\boA(t^4)$, whence
$\norm{\left[1+\frac{t^3}{1-tK(z)-t^3}\right]^n}_\boA \leq C
(1+2t^3)^n$. Let us estimate the coefficient of $z^p$ in $t^2
\sum_{n=0}^\infty (1-dt^2)^n
\left[1+\frac{t^3}{1-tK(z)-t^3}\right]^n z^n$. This is at most
  \begin{align*}
  t^2\sum_{n=0}^\infty (1-dt^2)^n
  \left(\left[1+\frac{t^3}{1-tK(z)-t^3}\right]^n\right)_{-n+p}
  &
  \leq \sum_{n=0}^\infty t^2(1-dt^2)^n
  \frac{(1+2t^3)^n}{|-n+p|^\gamma}
  \\&
  \leq \sum_{n=0}^\infty t^2\left(1-\frac{d}{2}t^2\right)^n
  \frac{1}{|-n+p|^\gamma}
  \end{align*}
for $t$ small enough so that $(1-dt^2)(1+2t^3)\leq
\left(1-\frac{d}{2}t^2\right)$. We find the same expression as in
the convolution between $O\left(1_{n\geq 0} t^2
\left(1-\frac{d}{2}t^2\right)^n\right)$ and $O(1/|n|^\gamma)$,
that we have already estimated in Lemma \ref{ngamma_star_t}. Thus,
we get at most $O\left(\frac{1}{|p|^\gamma} +1_{p\geq
0}t^2\left(1-\frac{d}{4}t^2\right)^p\right)$.

As $\frac{t z (1-dt^2)}{1-tK(z)-t^3} =O(t/|n|^\gamma)$, another
convolution yields that \eqref{aerj;sdf} is
  \begin{equation*}
  1+O\left(\frac{t}{|n|^\gamma} +1_{n\geq 0}
  t^3\left(1-\frac{d}{8}t^2\right)^n\right).
  \end{equation*}
Multiplying by
$\frac{1}{1-tK(z)-t^3}=1+O\left(\frac{t}{|n|^\gamma}\right)$ gives that
 \eqref{yosjadgh} $=1+O\left(\frac{t}{|n|^\gamma}
+1_{n\geq 0} t^3\left(1-\frac{d}{16}t^2\right)^n\right)$. This
concludes the proof in the case $F(z)=1$.

We now handle the case of an arbitrary $F(z)$. Note that
  \begin{equation*}
  \bigl[F(z)+G_t(z) \bigr]^{-1}
  =\bigl[1+F(z)^{-1}G_t(z)\bigr]^{-1} F(z)^{-1}.
  \end{equation*}

The Wiener Lemma \ref{wiener_S1} implies that
$F(z)^{-1}=O\left(\frac{1}{|n|^\gamma}\right)$, whence Lemma
\ref{ngamma_star_t} gives
$F(z)^{-1}G_t(z)=O\left(\frac{t}{|n|^\gamma}+1_{n\geq 0}
t^3\left(1-\frac{d}{2}t^2\right)^{n}\right)$. Thus, the case
$F(z)=1$ yields
  \begin{equation*}
  \bigl[1+F(z)^{-1}G_t(z)\bigr]^{-1}=1+O\left(\frac{t}{|n|^\gamma}+1_{n\geq
  0} t^3 \left(1-\frac{d}{32} t^2\right)^n\right).
  \end{equation*}
Another convolution with
$F(z)^{-1}$ gives the result.
\end{proof}

\subsection{Proof of Theorem \ref{thm_abstrait_vitesse}}

Let $M(t)$ be as in Theorem \ref{thm_abstrait_vitesse}. We fix
once and for all $d>0$ such that $|1-\frac{1}{\mu}M(t)| \leq
1-dt^2$ for small enough $t$, and we restrict the range of $t$ so
that this inequality is true. The invertibility of $R(z,t)$ for
$z\in \D$ and small enough $t$ is proved in \cite[Proposition
2.7]{gouezel:stable}.

To estimate the eigenvalues using Banach algebra techniques, we will
need that the eigenvalue close to $1$ of $R(z)$ is defined on the
whole circle $S^1$, which is not \emph{a priori }the case.
Consequently, we use the construction in the second step of the proof
of Theorem 2.4 in \cite{gouezel:stable}: we replace $R(z)$ by
$\tilde{R}(z)=\sum_{-\infty}^\infty \tilde{R}_n z^n$, such that it
has a unique eigenvalue $\tilde{\lambda}(z)$ close to $1$ for $z\in
S^1$, equal to $1$ only for $z=1$, with
$\sum_{|k|>n} \| \tilde{R}_k \|=O(1/n^{\beta})$,
and such that $R(z)=\tilde{R}(z)$
for $z$ close to $1$ on $S^1$. We also set
$\tilde{R}(z,t)=(\tilde{R}(z)-R(z))+\sum R_n(t)z^n$. For small enough
$t$, $\tilde{R}(z,t)$ has for all $z\in S^1$ a unique eigenvalue
$\tilde{\lambda}(z,t)$ close to $1$.

\begin{lem}
\label{DL_lambda}
We have
  \begin{equation*}
  \frac{1-\tilde{\lambda}(z,t)}{1-(1-\frac{1}{\mu}M(t))z}
  =\frac{1-\tilde{\lambda}(z)}{1-z}
  +O\left(\frac{t}{|n|^{\beta-1}}+1_{n\geq 0}t^3\left(1-\frac{d}{2}t^2\right)^n
  \right).
  \end{equation*}
\end{lem}
%On peut jouer sur les hypothèses, modifier l'hypothèse sur
%$\norm{R_n}$, pour obtenir des r\'esultats un peu diff\'erents
%dans les puissances de $n$ (et de $t$) qui apparaissent.
\begin{proof}
We write
$K(z,t)=\frac{\tilde{\lambda}(1,t)-\tilde{\lambda}(z,t)}{1-z} -
\frac{1-\tilde{\lambda}(z)}{1-z}$. Recall that
$\tilde{\lambda}(1,t)=\lambda(1,t)=1-M(t)$. Then, writing
$B(z)=\frac{1-\tilde{\lambda}(z)}{1-z}$, we have
  \begin{equation*}
  \tilde{\lambda}(z,t)=1-M(t)+(z-1)(K(z,t)+B(z)).
  \end{equation*}
Thus,
  \begin{align*}
  \frac{1-\tilde{\lambda}(z,t)}{1-(1-\frac{1}{\mu}M(t))z}-B(z)
  &=\frac{M(t)+(1-z)(K(z,t)+B(z)) -
  \bigl[1-(1-\frac{1}{\mu}M(t))z\bigr]B(z)}{1-(1-\frac{1}{\mu}M(t))z}\\
  &=K(z,t)
  \frac{1-z}{1-(1-\frac{1}{\mu}M(t))z}
  +M(t) \frac{1-z B(z)/\mu}{1-(1-\frac{1}{\mu}M(t))z}
  \\&
  =I+II.
  \end{align*}

For $I$,
  \begin{equation}
  \label{expansion_Mt}
  \frac{1-z}{1-(1-\frac{1}{\mu}M(t))z}
  =1-\frac{1}{\mu}M(t) \sum_{n=1}^\infty
  \left(1-\frac{1}{\mu}M(t) \right)^{n-1} z^n
  \end{equation}
is in $1+O(1_{n\geq 0} t^2 (1-d t^2)^n)$. We multiply it by
$K(z,t)$. Set $\boA=\boO_{\beta-1}(\Hom(E))$ (the Banach algebra
of functions whose coefficients are in $O(1/|n|^{\beta-1})$ for
$n\in \Z$). We have $\frac{\tilde{R}(z,t)-\tilde{R}(1,t)}{z-1}
=\frac{\tilde{R}(z)-\tilde{R}(1)}{z-1}+O_{\boA}(t)$. The proof of
\cite[Lemma 2.6]{gouezel:stable}, but in the algebra $\boA$ and
with tildes everywhere, applies (using Theorem \ref{wiener_S1} to ensure
that the inverses remain in $\boA$). It gives
$\frac{\tilde{\lambda}(z,t)
-\tilde{\lambda}(1,t)}{z-1}=\frac{\tilde{\lambda}(z)-
\tilde{\lambda}(1)}{z-1} +O_{\boA}(t)$, i.e.\
$K(z,t)=O(t/|n|^{\beta-1})$. Hence, Lemma
\ref{ngamma_star_t} yields
  \begin{equation*}
  I=O\left(\frac{t}{|n|^{\beta-1}}+1_{n\geq 0} t^3
  \left(1-\frac{d}{2}t^2\right)^{n}\right).
  \end{equation*}

Since $\frac{\tilde{R}(z)-\tilde{R}(1)}{z-1}=
O\left(\frac{1}{|n|^\beta}\right)$, we prove that $B(z)=
O\left(\frac{1}{|n|^\beta}\right)$ as in the third step of the
proof of Theorem 2.4 in \cite{gouezel:stable} (but in the Banach
algebra $\boO_\beta(\Hom(E))$). Since $1-zB(z) /\mu$ vanishes at $1$
(Step 7 of the proof of Lemma 3.1 in \cite{gouezel:decay}) and is
in $O\left(\frac{1}{|n|^\beta}\right)$, it can be written as
$(1-z)C(z)$ where $C(z)=O(1/|n|^{\beta-1})$.
\begin{comment}
Indeed, $1-zB(z)/\mu$ can be written
$\sum c_n (z^n-1)$ with $|c_n|\leq C/|n|^{\beta}$ (since
$\tilde{R}_n=O(1/|n|^{\beta+1})$). For $n\geq 0$,
$z^n-1=(z-1)\sum_{k=0}^{n-1}z^k$, and for $n<0$,
$z^n-1=(1-z)\sum_{k=n}^{-1}z^k$. Hence,
  \begin{equation*}
  1-\frac{zB(z)}{\mu}=(z-1)\left(-\sum_{n=-\infty}^{-1}
  \left(\sum_{k=-\infty}^n c_k\right)z^n + \sum_{n=0}^\infty
  \left(\sum_{k=n+1}^\infty c_k \right) z^n\right)
  \end{equation*}
and the coefficients on the righthand side are in
$O(1/|n|^{\beta-1})$.
\end{comment}
To obtain $II$, we multiply $C(z)$ by
$M(t)\frac{1-z}{1-(1-\frac{1}{\mu}M(t))z}=O(t^2+1_{n\geq
0}t^4(1-dt^2)^n)$.
Lemma \ref{ngamma_star_t} yields
  \begin{equation*}
  II=O\left(\frac{t^2}{|n|^{\beta-1}}+1_{n\geq 0} t^4
  \left(1-\frac{d}{2}t^2\right)^{n} \right).
  \qedhere
  \end{equation*}
\end{proof}

\begin{cor}
\label{estimees_tilde}
We have
  \begin{equation*}
  \left(\frac{I-\tilde{R}(z,t)}{1-\left(1-\frac{1}{\mu}M(t)\right)z}
  \right)^{-1}
  =\left(\frac{I-\tilde{R}(z)}{1-z}\right)^{-1}
  +O\left(\frac{t}{|n|^{\beta-1}}+1_{n\geq 0}
  t^3\left(1-\frac{d}{256}t^2\right)^n \right).
  \end{equation*}
\end{cor}
\begin{proof}
Let $\tilde{P}(z,t)$ be the spectral projection associated to the
eigenvalue $\tilde{\lambda}(z,t)$ of $\tilde{R}(z,t)$, and
$\tilde{Q}(z,t)=I-\tilde{P}(z,t)$. Set $\boA=\boO_{\beta-1}(\Hom(E))$.
Then $\tilde{R}(z,t)=\tilde{R}(z)+O_\boA(t)$ by assumption. As
$\boA$ satisfies the Wiener Lemma \ref{wiener_S1}, the integral
expression of the projection $\tilde{P}(z,t)$ shows that
$\tilde{P}(z,t)=\tilde{P}(z) +O_\boA(t)$. Moreover,
$I-\tilde{R}(z,t)\tilde{Q}(z,t)=I-\tilde{R}(z)\tilde{Q}(z)+O_\boA(t)$,
whence
$(I-\tilde{R}(z,t)\tilde{Q}(z,t))^{-1}=(I-\tilde{R}(z)\tilde{Q}(z))^{-1}
+O_\boA(t)$, since $I-\tilde{R}(z) \tilde{Q}(z)$ is invertible in
$\boA$ by Theorem \ref{wiener_S1}, and the inversion is Lipschitz.

As
  \begin{equation*}
  I-\tilde{R}(z,t)=(1-\tilde{\lambda}(z,t))\tilde{P}(z,t)
  +(I-\tilde{R}(z,t)\tilde{Q}(z,t)) \tilde{Q}(z,t),
  \end{equation*}
Lemmas \ref{DL_inversion} and \ref{DL_lambda} yield that
  \begin{align*}
  \left(\frac{I-\tilde{R}(z,t)}{1-(1-\frac{1}{\mu}M(t))z}\right)^{-1}
  \!\!\!\!\!\!\!\!\!\!\!\!\!\!\!\!\!\!\!\!\!\!\!\!\!\!\!\!\!\!
  \!\!\!\!\!\!\!\!\!\!\!\!\!\!\!\!
  \!\!\!\!\!\!\!\!\!\!\!\!\!\!&
  \\&
  =\frac{1-(1-\frac{1}{\mu}M(t))z}{1-\tilde{\lambda}(z,t)}\tilde{P}(z,t)
  +\left(1-\left(1-\frac{1}{\mu}M(t)\right)z\right)
  (I-\tilde{R}(z,t)\tilde{Q}(z,t))^{-1}
  \tilde{Q}(z,t)
  \\&
  =\left[\frac{1-z}{1-\tilde{\lambda}(z)}
  +O\left(\frac{t}{|n|^{\beta-1}}+1_{n\geq 0} t^3
  \left(1-\frac{d}{128} t^2\right)^n \right) \right]
  \left[ \tilde{P}(z)+O\left(\frac{t}{|n|^{\beta-1}}\right) \right]
  \\&\hphantom{=\ }
  +\left[ (1-z)+O(1_{n=0}t^2)\right]
  \left[ (I-\tilde{R}(z)\tilde{Q}(z))^{-1}
  +O\left(\frac{t}{|n|^{\beta-1}}\right) \right]
  \left[ \tilde{Q}(z)+O\left(\frac{t}{|n|^{\beta-1}}\right) \right]
  \\
  &=\frac{1-z}{1-\tilde{\lambda}(z)}\tilde{P}(z)
  + (1-z) (I-\tilde{R}(z)\tilde{Q}(z))^{-1} \tilde{Q}(z)
  +O\left( \frac{t}{|n|^{\beta-1}} +1_{n\geq 0}t^3 \left(1-\frac{d}{256}t^2
  \right)^n \right)
  \\
  &= \left(\frac{I-\tilde{R}(z)}{1-z}\right)^{-1}
  +O\left( \frac{t}{|n|^{\beta-1}} +1_{n\geq 0}t^3 \left(1-\frac{d}{256}t^2
  \right)^n \right).
  \qedhere
  \end{align*}
\end{proof}

\begin{cor}
\label{cor36}
We have
  \begin{equation*}
  \left(\frac{I-R(z,t)}{1-\left(1-\frac{1}{\mu}M(t)\right)z}
  \right)^{-1}
  =\left(\frac{I-R(z)}{1-z}\right)^{-1}
  +O\left(\frac{t}{|n|^{\beta-1}}+1_{n\geq 0}
  t^3\left(1-\frac{d}{512}t^2\right)^n \right).
  \end{equation*}
\end{cor}
\begin{proof}
Let $\chi_1,\chi_2$ be a $C^\infty$ partition of unity of $S^1$
such that $R(z)=\tilde{R}(z)$ on the support of $\chi_1$.
Since $\chi_1$ is $C^\infty$, $\chi_1(z)=O(1/|n|^{\beta-1})$.
Writing
$A(z,t)=\left(\frac{I-R(z,t)}{1-\left(1-\frac{1}{\mu}M(t)
\right)z}\right)^{-1}$, Corollary \ref{estimees_tilde} ensures
that
  \begin{equation*}
  \chi_1(z)A(z,t)=\chi_1(z)A(z,0)
  +O\left(\frac{t}{|n|^{\beta-1}}+1_{n\geq 0}
  t^3\left(1-\frac{d}{512}t^2\right)^n \right).
  \end{equation*}
Concerning $\chi_2$, we can modify $R$ outside of its support so that
$I-R(z)$ is everywhere invertible on $S^1$. Let
$\boA=\boO_{\beta-1}(\Hom(E))$.
Since $I-R(z,t)=I-R(z)+O_\boA(t)$, Theorem
\ref{wiener_S1} yields
$(I-R(z,t))^{-1}=(I-R(z))^{-1}+O_\boA(t)$. Hence,
  \begin{equation*}
  \chi_2(z)A(z,t)=\chi_2(z) A(z,0)+O\left(\frac{t}{|n|^{\beta-1}}\right).
  \qedhere
  \end{equation*}
\end{proof}

In fact, the functions in the previous Corollary are defined on
the whole disk $\overline{\D}$, i.e.\ their coefficients for $n<0$
vanish. However, during the proof, it was important to work in a
less restrictive context, for example to introduce partitions of
unity.

\begin{proof}[Proof of Theorem \ref{thm_abstrait_vitesse}]
Set
  \begin{equation*}
  F_t(z)=\left(\frac{I-R(z,t)}{1-\left(1-\frac{1}{\mu}M(t)\right)z}
  \right)^{-1}
  -\left(\frac{I-R(z)}{1-z}\right)^{-1}
  \end{equation*}
and
  \begin{equation*}
  \left(\frac{I-R(z)}{1-z}\right)^{-1}=\frac{1}{\mu}P+(1-z)A(z)
  \end{equation*}
where $A(z)=O\left(\frac{1}{n^{\beta-1}}\right)$, by \cite[Theorem
1]{sarig:decay} or \cite[Theorem 1.1]{gouezel:decay}.

Then
  \begin{align*}
  \sum T_{n}(t)z^n&=(I-R(z,t))^{-1}
  \\&
  =\frac{1}{1-(1-\frac{1}{\mu}M(t))z} \frac{1}{\mu}P
  +\frac{1-z}{1-(1-\frac{1}{\mu}M(t))z}A(z)
  +\frac{1}{1-(1-\frac{1}{\mu}M(t))z}F_t(z)
  \\&
  =I+II+III.
  \end{align*}
The coefficient of $z^n$ in $I$ is
$\frac{1}{\mu}\left(1-\frac{1}{\mu}M(t) \right)^n P$. So, we have
to bound the coefficients of $II$ and $III$ to conclude.

By \eqref{expansion_Mt} and Lemma \ref{ngamma_star_t}, the
coefficients of $II$ belong to
  \begin{equation*}
  \left[ 1+O\left( t^2 (1-dt^2)^n \right) \right]\star
  O\left( \frac{1}{n^{\beta-1}}\right)
  \subset
  O\left(\frac{1}{n^{\beta-1}}+t^2\left(1-\frac{d}{2}t^2\right)^n
  \right).
  \end{equation*}

For $III$, we get by Corollary \ref{cor36} that the coefficient of
$z^n$ is bounded by
  \begin{equation*}
  O\left( (1-dt^2)^n \right) \star O\left(\frac{t}{n^{\beta-1}}+
  t^3\left(1-\frac{d}{512}t^2\right)^n \right).
  \end{equation*}

The convolution between $(1-dt^2)^n$ and $t^3
(1-\frac{d}{512}t^2)^n$ is bounded by the convolution between $
(1-\frac{d}{512}t^2)^n$ and $t^3
 (1-\frac{d}{512}t^2)^n$, which gives $n t^3
(1-\frac{d}{512}t^2)^n$. This is less than $C t
(1-\frac{d}{1024}t^2)^n$, since
  \begin{equation*}
  \frac{n}{2} (1-ct^2)^{n} \leq \sum_{k=0}^{n/2} (1-ct^2)^k
  (1-ct^2)^{n/2} \leq (1-ct^2)^{n/2} \frac{1}{ct^2}.
  \end{equation*}
As $t (1-\frac{d}{1024}t^2)^n \leq t
(1-\frac{d}{1024}t^2)^n \star \left(\frac{1}{n^{\beta-1}}\right)$,
we get a bound
of the form stated in the theorem.
\end{proof}

\section{The key estimate}
\label{section_estimees_tour}
In this section, the assumptions are as in Theorem
\ref{TCL_normal}, i.e.\ $X$ is a Young tower whose return time
$\phi$ satisfies $m[\phi>n]=O(1/n^{\beta})$ with $\beta>2$. We
also assume that $\gcd(\phi_i)=1$. Take also $f\in C_\tau(X)$ with
$\int f \dd m=0$.

The goal of this section is to estimate precisely $\int_X e^{it
S_n f} \cdot u \cdot v\circ T^n \dd m$
when $u,v$ are functions on $X$. We will use the
same kind of perturbative ideas as in \cite{rousseau-egele} and
\cite{guivarch-hardy}, but applied to transfer operators
associated to induced maps on the basis $B$ of the tower
(see \cite{sarig:decay}). Separating the different return times, it
will be possible to use the abstract Theorem
\ref{thm_abstrait_vitesse}, to finally get the key estimate
Theorem \ref{estimee_clef_markov}.

\subsection{First return transfer operators}

Let $\hat{T}$ be the transfer operator associated to $T$, defined
by $\int u\cdot v\circ T \dd m=\int \hat{T}u \cdot v \dd m$. When
$\psitest$ is integrable on $X$, it is given by
  \begin{equation*}
  \hat{T}^n \psitest(x)=\sum_{T^n y=x} g_m^{(n)}(y) \psitest(y),
  \end{equation*}
where $g_m^{(n)}$ is the inverse of the jacobian of $T^n$ at $y$.

As the basis $B$ of the tower plays a particular role, we will
decompose the trajectories of the preimages of $x$ under $T^n$,
keeping track of the moments their iterates fall again in $B$.
More formally, we introduce the following operators:
  \begin{gather*}
  R_n \psitest(x)=\sum_{\substack{T^n y=x\\
    y\in B, Ty,\ldots,T^{n-1}y \not \in B, T^n
    y\in B}} g_m^{(n)}(y) \psitest(y),
  \\
  T_n\psitest(x)=\sum_{\substack{T^n y=x\\
    y\in B, T^n y\in B}}
    g_m^{(n)}(y)\psitest(y),
  \\
  A_n \psitest(x)=\sum_{\substack{T^n y=x\\
   y\in B, Ty,\ldots,T^n y\not\in B}}
   g_m^{(n)}(y)\psitest(y),
  \\
  B_n \psitest(x)=\sum_{\substack{T^n y=x\\
   y,\ldots,T^{n-1}y\not \in B, T^n y\in B}}
   g_m^{(n)}(y)\psitest(y),
  \\
   C_n \psitest(x)=\sum_{\substack{T^n y=x\\
   y,\ldots,T^{n}y\not \in B}}
   g_m^{(n)}(y)\psitest(y).
  \end{gather*}
The interpretations of these operators are as follows: $R_n$ takes
into account the first returns to $B$, while $T_n$ takes all
returns into account. Hence,
  \begin{equation}
  \label{renouvellement}
  T_n=\sum_{k_1+\ldots+k_l=n} R_{k_1}\dots R_{k_l}.
  \end{equation}
The operators $B_n$ and $A_n$ see respectively the beginning and
the end of the trajectories, outside of $B$. Thus, if $x$ is fixed
and $y$ satisfies $T^n y=x$ and $\{y,Ty,\dots, T^n y\}\cap B \not=
\emptyset$, we can consider the first $b$
iterates of $y$, until it enters in $B$ (this corresponds to
$B_b$), then some successive returns to $B$, during a time $k$
(this corresponds to $T_k$), and finally $a$ iterations outside of
$B$ (corresponding to $A_a$). Thus,
  \begin{equation}
  \label{exprime_operateur_transfert}
  \hat{T}^n=C_n+\sum_{a+k+b=n} A_a T_k B_b.
  \end{equation}
The operator $C_n$ takes into account the points $y$ such that
$\{y,Ty,\dots, T^n y\}\cap B=\emptyset$.

We perturb these operators, setting (for $X=A,B,C,R,T$, and $t\in
\R$)
  \begin{equation*}
  X_n(t)(\cdot)=X_n(e^{itS_n f}\cdot).
  \end{equation*}
Equation \eqref{renouvellement} remains true with $t$
everywhere:
  \begin{equation}
  \label{renouvellement_perturbe}
  T_n(t)=\sum_{k_1+\ldots+k_l=n} R_{k_1}(t)\dots R_{k_l}(t).
  \end{equation}
Let $\hat{T}(t)$ be the perturbation of $\hat{T}$ given by
$\hat{T}(t)(\cdot)=\hat{T}(e^{itf} \cdot)$. Then the following analogue
of \eqref{exprime_operateur_transfert} holds:
  \begin{equation}
  \label{exprime_operateur_transfert_perturbe}
  \hat{T}(t)^n=C_n(t)+\sum_{a+k+b=n} A_a(t) T_k(t) B_b(t).
  \end{equation}

Let $f_B$ be given by \eqref{definition_fY}. As $f_B=S_n f$ on
$\{y\in B \tq \phi(y)=n\}$, we get $R_n(t) (\psitest)=R_n(e^{it f_B}
\psitest)$.

For $z\in \overline{\D}$, write $R(z)=\sum R_n z^n$, and
  \begin{equation}
  \label{definit_Rzt}
  R(z,t)=\sum_{n=0}^\infty R_n(t) z^n.
  \end{equation}
Let $T_B$ be the first return map induced by $T$ on $B$, i.e.\
$T_B(x)=T^{\phi(x)}(x)$. Then \eqref{definit_Rzt}
corresponds to considering all the preimages of a point in
$B$ under $T_B$, whence $R(1)$ is the transfer operator
$\hat{T}_B$ associated to $T_B$, and
  \begin{equation}
  \label{operateur_transfert_induit}
  R(1,t)(\psitest)=
  \hat{T}_B (e^{it f_B} \psitest).
  \end{equation}

When $\psitest$ is a function on a subset $Z$ of $X$, we denote by
$D_\tau \psitest(Z)$ the best H\"older constant of $\psitest$ on $Z$, i.e.\
  \begin{equation}
  \label{definit_Dtau}
  D_\tau \psitest(Z)=\inf \{C>0 \tq \forall x,y \in Z, |\psitest(x)-\psitest(y)|\leq C
  \tau^{s(x,y)} \}
  \end{equation}
where $s(x,y)$ is the separation time of $x$ and $y$.

The operators $R_n$ and $R_n(t)$ act on $C_{\tau'}(B)$ for any
$\tau\leq \tau' <1$. Take $\eta$ such that
  \begin{equation}
  \label{definit_eta}
  0<\eta<\min(1/2, \beta-2)
  \end{equation}
and set
$\nu=\tau^\eta$.
For technical reasons, we will let the operators act on
$C_{\nu}(B)$. We regroup in
the following lemma all the estimates we will need later. Their
proofs are rather straightforward, but sometimes lengthy. Hence, we
will not give the details of the proofs, and rather
give references to articles where similar estimates
are proved.
\begin{lem}
\label{enumere_estimees}
The operators $R_n$ and $R_n(t)$ acting on $C_{\nu}(B)$
satisfy the following estimates:
\begin{enumerate}
\item $\sum_{k=n+1}^\infty \norm{R_k}=O(1/n^{\beta})$.
\item
\label{point_2_lemme}
The operator $R(z,t)$ satisfies a Doeblin-Fortet inequality
  \begin{equation}
  \label{doeblin}
  \norm{R(z,t)^n u} \leq C(1+|t|)(\nu |z|)^n \norm{u}+C|z|^n
  \norm{u}_{L^1}.
  \end{equation}
In particular, the spectral radius of $R(z,t)$ is $\leq |z|$, and its essential
spectral radius is $\leq \nu |z|$. Thus, $I-R(z,t)$ is not
invertible if and only if $1$ is an eigenvalue of $R(z,t)$, and
this can happen only for $|z|=1$.
\item
\label{point_3_lemme}
$R(1)$ has a simple isolated eigenvalue at $1$ (the eigenspace is
the space of constant functions), and $I-R(z)$ is invertible for $z \not
=1$.
\item
\label{point_4_lemme}
There exists $C>0$ such that, for any $t\in \R$, for every
$n\in \N^*$,
  \begin{equation*}
  \sum_{k=n+1}^\infty
  \norm{R_k(t)-R_k(0)} \leq C\frac{ |t|}{n^{\beta-1}}.
  \end{equation*}
\item
\label{point_5_lemme}
For every $t\in \R$, there exists $C=C(t)$ such that, for any $t'\in
\R$, for every $n\in \N^*$, $\norm{R_n(t')-R_n(t)} \leq C
\frac{|t'-t|^{1/2}}{n^{\beta-1}}$.
\end{enumerate}
\end{lem}
\begin{proof}
The first assertion is a consequence of \cite[Lemma
4.2]{gouezel:stable}: it gives $\norm{R_n}=O(m[\phi=n])$.

The inequality \eqref{doeblin} is similar to \cite[Proposition
2.1]{aaronson_denker}. In this article, the hypothesis is that
$D_\nu f_B(B)<\infty$, but the proofs work in fact as soon as $\sum
m(B_{i,0}) D_\nu f_B(B_{i,0})<\infty$. In our case, $D_\nu f_B(
B_{i,0}) \leq C \phi_i$, whence $\sum m(B_{i,0}) D_\nu
f_B(B_{i,0}) \leq C \sum m(B_{i,0}) \phi(B_{i,0})=C$ by Kac's
Formula. Since the injection $C_{\nu}(B) \to L^1(B)$ is
compact, the statement on the essential spectral radius of $R(z,t)$ is
then a consequence of Hennion's Theorem (\cite{hennion}).

The third assertion is \cite[Lemma 4.3]{gouezel:stable}.

The two remaining assertions are proved by direct estimates, similar
to the estimates in
\cite[Theorem 2.4]{aaronson_denker}. The last one holds in
$C_{\sqrt{\tau}}(X)$ but not in $C_\tau(X)$.
This is the reason of the requirement $\nu\geq \sqrt{\tau}$.
\end{proof}

We will also need the following estimates on $A_a(t)$,
$B_b(t)$ and $C_n(t)$, acting on $C_{\nu}(X)$.

\begin{lem}
\label{int_Aa}
Let $u,v\in L^\infty(X)$. There exists a constant $C$ such that, for
any $t\in \R$, for any $a\in \N$,
  \begin{equation}
  \label{int_Aa1}
  \left|\int A_a(t)(u)v -  \int A_a(u)v\right|\leq
C\frac{|t|}{a^{\beta-1}}\norm{u}_\infty \norm{v}_\infty.
  \end{equation}
Moreover, $\int A_a(u)v =O(1/a^\beta)$, and
  \begin{equation}
  \label{int_Aa2}
  \sum_{a=0}^\infty \int A_a(1_B) v = \int v.
  \end{equation}
\end{lem}
\begin{proof}
Define a function $v'$ on $B$ by $v'(x)=v(T^a x)$ if $\phi(x)>a$, and
$0$ otherwise.
Changing variables, we get $\int A_a(t)(u)v=\int_{\{\phi>a\}}
e^{it S_a f}u v'$. Since $|e^{it S_a f}-1|\leq |t| a \norm{f}_\infty$
and $m( \phi>a) \leq \frac{C}{a^\beta}$, this implies \eqref{int_Aa1}.

Moreover, $\left| \int A_a(u)v \right| \leq m(\phi>a) \norm{u}_\infty
\norm{v}_\infty \leq \frac{C}{a^\beta} \norm{u}_\infty
\norm{v}_\infty$. Finally, $ \int A_a(1_B) v = \int_{T^a\{\phi>a\}}
v$. Since the sets $T^a\{\phi>a\}$ form a partition of $X$,
\eqref{int_Aa2} readily follows.
\end{proof}

\begin{lem}
\label{int_Bb} For $t\in \R$, $B_b(t)=B_b +O\left(
\frac{|t|}{b^{\beta-1}}\right)$, where $\norm{B_b}=O(1/b^{\beta})$
and $\forall \psitest \in C_{\nu}(X)$,
  \begin{equation*}
  \sum_{b=0}^{\infty} \int_B B_b \psitest\dd m =\int_X \psitest \dd m.
  \end{equation*}
\end{lem}
\begin{proof}
Let $\Lambda_b$ be the set of points that enter into $B$ after
exactly $b$ iterations, so that $B_b(t)(\psitest)$ takes the values of
$\psitest$ on $\Lambda_b$ into account. As in \cite[Theorem
2.4]{aaronson_denker}, we check that $\norm{B_b(t)-B_b} \leq C|t|
m(\Lambda_b) D_\nu S_b f(\Lambda_b)$. As $D_\nu S_b f\leq
C b$ and
$m(\Lambda_b)=m[\phi>b]=O(1/b^\beta)$, we get indeed that
$\norm{B_b(t) -B_b} \leq C \frac{|t|}{b^{\beta-1}}$.

Moreover, we check as in  \cite[Lemma 4.2.]{gouezel:stable} that
$\norm{B_b}=O(m(\Lambda_b))=O(1/b^\beta)$. Finally, as $\int_B B_b
\psitest \dd m =\int_{\Lambda_b} \psitest \dd m$ we get $\sum_b \int_B B_b
\psitest\dd m = \int_X \psitest \dd m$.
\end{proof}

\begin{lem}
\label{int_Cn}
Let $u\in L^\infty(X)$. There exists a constant $C>0$ such that, for
any $t\in \R$, for any $n\in \N$,
  \begin{equation*}
  \norm{C_n(t)(u)}_{L^1} \leq \frac{C}{n^{\beta-1}} \norm{u}_\infty.
  \end{equation*}
\end{lem}
\begin{proof}
Set $X_{n+1}= X \moins \bigcup_{i=0}^n T^{-i}B$ and
$Z_{n+1}=T^n(X_{n+1})$. The function $C_n(u)$ vanishes outside of
$Z_{n+1}$. Let $x\in Z_{n+1}$ and let $y$ be its preimage under $T^n$.
Then
$|C_n(t)(u)(x)|=|u(y)| \leq \norm{u}_\infty$. Hence,
$\norm{C_n(t)(u)}_{L^1} \leq m(Z_{n+1}) \norm{u}_\infty$.
Finally, $m(Z_{n+1})=m(X_{n+1})=\sum_{p=n+1}^\infty m(
\phi>p)=O(1/n^{\beta-1})$.
\end{proof}

\subsection{Result for the induced map}

Equation \eqref{operateur_transfert_induit} and
\ref{point_2_lemme} in Lemma \ref{enumere_estimees} imply that the
perturbed transfer operator $R(1,t)=\hat{T}_B(t)=\hat{T}_B(e^{it
f_B} \cdot)$ has a spectral gap. Thus, the classical methods of
\cite{guivarch-hardy} apply
to it, and yield an asymptotic expansion of the maximal eigenvalue
of $\hat{T}_B(t)$ (this implies a central limit theorem for
$T_B$ and the function $f_B$,
but we are not interested in it here). To estimate the speed in
the central limit theorem, we will need a rather precise asymptotic
expansion of this eigenvalue, given by the following proposition.

\begin{prop}
\label{TCL_induit}
Assume that $T$ and $f$ satisfy the assumptions of Theorem
\ref{TCL_normal} in a mixing Young tower.
Let $\sigma^2$ be the variance in this theorem.
Then, for small enough $t$, $R(1,t)$ has a unique eigenvalue
$\lambda(1,t)$ close to $1$. It can be written as
$\lambda(1,t)=1-\frac{\sigma^2}{2m(B)} L(t)$ for a function $L$ such
that $L(t)\sim t^2$ when $t\to 0$.

Write $E_B(u)=\frac{1}{m(B)}\int_B u$,
and define a function $a$ on $B$ by $a=(I-\hat{T}_B)^{-1}( \hat{T}_B f_B)$.
Then
  \begin{equation}
  \label{dev}
  \lambda(1,t)=E_B(e^{it f_B}) - t^2 E_B (a f_B) +O(t^3).
  \end{equation}
\end{prop}
\begin{proof}
The fact that $\lambda(1,t)=1-\frac{\sigma^2}{2m(B)}(t^2+o(t^2))$ is a
consequence of
\cite[Theorem 3.7]{gouezel:stable}. It remains to prove \eqref{dev}.
As everything takes place in $B$, we can
multiply $m$ by a constant, and assume that $m(B)=1$. Set
$R_t=R(1,t)$, and let $\xi_t$ be the eigenfunction of $R_t$
corresponding to its eigenvalue $\lambda_t$ close to $1$. We
normalize it so that $\int \xi_t=1$. We will also write $R=R_0=\hat{T}_B$.

Lemma 3.4 of
\cite{gouezel:stable} states that there exists a constant $C$ such
that, if $g:B \to \R$ is integrable,
  \begin{equation}
  \label{lemme34_ancien}
  \norm{R g}_\nu \leq C \left(\norm{g}_{L^1}+\sum_I m(B_{i,0}) D_\nu
  g(B_{i,0}) \right)
  \end{equation}
where $D_\nu g(B_{i,0})$ is the best $\nu$-H\"older constant of
$g$ on $B_{i,0}$,
defined in
\eqref{definit_Dtau}. In particular, $Rf_B\in C_\nu(B)$.

Let us show that, in $C_\nu(B)$,
  \begin{equation}
  \label{equation_dans_Cnu}
  R \left(\frac{e^{itf_B}-1}{t} \right) = iR(f_B) +O(t).
  \end{equation}
The Taylor Formula gives
  \begin{equation*}
  \frac{e^{it f_B}-1}{t} -i f_B = -t f_B^2 \int_0^ 1
  (1-u) e^{it u f_B} \dd u.
  \end{equation*}
We use this formula to bound $D_\nu\left(\frac{e^{it f_B}-1}{t} -i
f_B \right) (B_{i,0})$, where $B_{i,0}$ is an element of the
partition of $B$. Set $n=\phi(B_{i,0})$ the return time on
$B_{i,0}$. As $D_\nu(h_1 h_2)\leq D_\nu(h_1) \norm{h_2}_\infty +
\norm{h_1}_\infty D_\nu(h_2)$,
  \begin{align*}
  D_\nu\left(\frac{e^{it
  f_B}-1}{t} -i f_B \right) (B_{i,0})
  & \leq
  |t| \norm{(f_B)_{|B_{i,0}}}_\infty D_\nu f_B (B_{i,0})
  \\& \hphantom{\leq \ }
  + |t| \norm{(f_B)_{|B_{i,0}}}_\infty^2 \int_0^1 D_\nu(e^{it u
  f_B}) (B_{i,0}) \dd u.
  \end{align*}
The first term is $\leq |t|n^2$. For the second term, take $C$ such
that  $|e^{is}-1| \leq C |s|^\eta$ for any $s\in \R$ (where $\eta$ was
defined in \eqref{definit_eta}).
\begin{comment}
(c'est facile à garantir pour $s$ petit ou grand, donc c'est
possible par compacit\'e).
\end{comment}
Then, if $x,y\in B_{i,0}$,
  \begin{align*}
  \left|e^{itu f_B(x)} -e^{itu f_B(y)} \right|&
  =\left| e^{itu (f_B(x)-f_B(y))}-1 \right|
  \leq C|tu (f_B(x)-f_B(y))|^\eta
  \\&
  \leq C | D_\tau f_B(B_{i,0}) \tau^{s(x,y)} |^\eta
  \leq C n^\eta \nu^{s(x,y)},
  \end{align*}
whence $D_\nu(e^{it u
  f_B})(B_{i,0}) \leq C n^\eta$.
An integration yields
  \begin{equation*}
  D_\nu\left(\frac{e^{it
  f_B}-1}{t} -i f_B \right) (B_{i,0})
  \leq |t| n^{2+\eta}.
  \end{equation*}
Equation \eqref{lemme34_ancien} gives that
  \begin{equation*}
  \norm{ R\left(\frac{e^{it
  f_B}-1}{t} -i f_B \right)}_\nu
  \leq C \left( \norm{\frac{e^{it
  f_B}-1}{t} -i f_B}_{L^1} + \sum m[\phi =n] |t| n^{2+\eta}\right).
  \end{equation*}
As $|e^{it f_B}-1 -it f_B| \leq t^2 f_B^2$, with $f_B^2$
integrable, the first term is $O(t)$. For the second term,
$ \sum m[\phi =n] n^{2+\eta} = \sum \bigl( m[\phi>n-1]- m[\phi>n]\bigr)
n^{2+\eta} \leq C \sum m[\phi>n] n^{1+\eta}$. Since
$m[\phi>n]=O(1/n^\beta)$, this sum is finite by
definition of $\eta$. This proves \eqref{equation_dans_Cnu}.

We return to the study of the eigenvalue $\lambda_t$.
As $\lambda_t \xi_t=R_t \xi_t$, we get after
integration that
  \begin{equation}
  \label{developpement_lambdat}
  \lambda_t = E(e^{it f_B}) +\int (e^{it f_B}-1)(\xi_t -1).
  \end{equation}
As $f_B \in L^2$ and $\int f_B=0$, we have
  \begin{equation}
  \label{developpement_dans_L3}
  E(e^{it f_B})=1+it \int f_B -\frac{t^2}{2} \int f_B^2+o(t^2)
  =1 -\frac{t^2}{2} \int f_B^2+o(t^2).
  \end{equation}
Moreover, by \ref{point_4_lemme} in Lemma
\ref{enumere_estimees}, $\xi_t-1=O(\norm{R_t-R})=O(t)$ in
$C_\nu(B)$, hence in $L^2$, and $e^{itf_B}-1=it f_B+o(t)=O(t)$
in $L^2$. Consequently, $\int (e^{it f_B}-1)(\xi_t -1)=O(t^2)$, which
implies that
$\lambda_t=1+O(t^2)$.

Thus,
  \begin{equation*}
  \frac{\xi_t-1}{t}=\frac{\lambda_t\xi_t -\xi_0}{t}+O(t)
  =(R_t-R_0)\frac{\xi_t-\xi_0}{t} +R \frac{\xi_t-\xi_0}{t}
  +\frac{R_t-R_0}{t} \xi_0 +O(t).
  \end{equation*}
As $R_t-R_0=O(t)$ and $\frac{\xi_t-\xi_0}{t}$ is bounded,
$(R_t-R_0)\frac{\xi_t-\xi_0}{t}=O(t)$. Moreover,
$\frac{R_t-R_0}{t} \xi_0 = R\left( \frac{e^{itf_B}-1}{t}\right)$.
Hence,
  \begin{equation*}
  (I-R) \frac{\xi_t -\xi_0}{t} = R\left(
  \frac{e^{itf_B}-1}{t}\right)+O(t).
  \end{equation*}
As $Rf_B \in C_\nu(B)$, and $I-R$ is invertible on the functions
of $C_\nu(B)$ with vanishing integral (\ref{point_3_lemme} in
Lemma \ref{enumere_estimees}), it is possible to define
$a=(I-R)^{-1} (R f_B) \in C_\nu(B)$. Then, using
\eqref{equation_dans_Cnu},
  \begin{equation*}
  (I-R) \left( \frac{\xi_t -\xi_0}{t} -ia \right)
  =R \left(\frac{e^{itf_B}-1}{t} -i f_B \right) +O(t)
  =O(t).
  \end{equation*}
As the inverse of $I-R$ is continuous on the functions with zero
integral, we get that, in $C_\nu(B)$ (hence in $L^\infty$),
  \begin{equation*}
  \xi_t=1+tia+O(t^2).
  \end{equation*}
As $e^{itf_B}=1+itf_B+O(t^2)$ in $L^1$ since $f_B^2$ is
integrable, we get
  \begin{equation*}
  \int (e^{itf_B}-1)(\xi_t-1) =-t^2 \int f_B a +O(t^3).
  \end{equation*}
Equation \eqref{developpement_lambdat} yields the desired conclusion.
\end{proof}

\subsection{The key estimate}

\begin{thm}
\label{estimee_clef_markov}
Let $X$ be a Young tower with $m[\phi>n]=O(1/n^{\beta})$ for
$\beta>2$, and $\gcd(\phi_i)=1$. Let $f\in C_\tau(X)$ be of zero
integral. Let $\sigma$ and $L(t)$ be given by Proposition
\ref{TCL_induit}, and assume $\sigma>0$. Then there exist
$\alpha>0$, $C>0$ and $d>0$ such that, for any $u \in C_\tau(X)$ and
$v\in L^\infty(X)$,
for any $n\in \N^*$, for any $t \in [-\alpha,\alpha]$,
  \begin{multline*}
  \left|\int_X e^{it S_n f} \cdot u\cdot v\circ T^n \dd m -
  \left( 1-\frac{\sigma^2}{2} L(t)
  \right)^n \left(\int_X u\dd m\right) \left(\int_X v\dd m\right)
  \right|
  \\
  \leq C \left[\frac{1}{n^{\beta-1}}+|t|\left(\frac{1}{n^{\beta-1}}\right)
  \star (1-dt^2)^n \right] \norm{u} \norm{v}_\infty.
  \end{multline*}
\end{thm}
\begin{proof}
Set $R(z,t)=\sum R_n(t) z^n$, we want to apply Theorem
\ref{thm_abstrait_vitesse} to $R(z,t)$. Proposition
\ref{TCL_induit} gives the behavior of the eigenvalue of $R(1,t)$,
while Lemma \ref{enumere_estimees} shows the required estimates.
Finally, the spectral projection $P$ of $R(1,0)$ is the projection
on the constant functions on $B$, given by $Pg=\frac{\int_B g \dd
m}{m(B)}$. It satisfies $PR'(1)P=\frac{1}{m(B)}P$ (\cite[Lemma
4.4.]{gouezel:stable}).

Consequently, we can apply
Theorem \ref{thm_abstrait_vitesse} with
$M(t)=\frac{\sigma^2}{2m(B)}L(t)$ and $\mu=\frac{1}{m(B)}$.
As $\sum T_n(t)
z^n =(I-R(z,t))^{-1}$ by \eqref{renouvellement_perturbe}, we get
that there exists an error term $E(n,t)$ such that $T_n(t)=m(B)
\left(1- \frac{\sigma^2}{2}L(t)\right)^n P+E(n,t)$, with
  \begin{equation}
  \norm{E(n,t)} \leq
  C\left[\frac{1}{n^{\beta-1}}
  +|t| \left( \frac{1}{n^{\beta-1}} \right) \star
  \left(1-dt^2\right)^n \right]=: e(n,t).
  \end{equation}
In the following, $C$ and $d$ will denote generic constants, that
may vary finitely many times. In particular, we may write
inequalities like $10 e(n,t) \leq e(n,t)$. With this convention,
Lemma \ref{ngamma_star_t} implies that
  \begin{equation}
  \label{convoles_ent}
  \left(\frac{1}{n^{\beta-1}} \right) \star e(n,t) = O(e(n,t)).
  \end{equation}

Let $u\in C_\tau(X)$ and $v\in L^\infty(X)$.
To simplify the expressions, we will assume that these functions are
of norm at most $1$.
Then
\eqref{exprime_operateur_transfert_perturbe} implies that
  \begin{align}
  \label{decompose_somme_avec_Ent}
    \begin{split}
  \int_X e^{it S_n f}\cdot u \cdot v\circ T^n &
  =\int_X \hat{T}^n (e^{it S_n f} u) v
  =\int_X \hat{T}(t)^n (u)v
  \\&
  =
  \sum_{a+k+b=n} \int_X A_a(t) T_k(t) B_b(t) (u)v + \int_X C_n(t)(u)v
  \\ &
  =\sum_{a+k+b=n} \int_X A_a(t) m(B)
  \left(1- \frac{\sigma^2}{2}L(t)\right)^k P B_b(t) (u)v
  \\& \hphantom{=\ }
  +\sum_{a+k+b=n} \int_X A_a(t) E(k,t) B_b(t) (u)v+\int_X C_n(t)(u)v.
  \end{split}
  \end{align}
By Lemma \ref{int_Cn}, $\left| \int_X C_n(t)(u)v \right| \leq
\frac{C}{n^{\beta-1}}\leq e(n,t)$.

Let us bound the second sum of \eqref{decompose_somme_avec_Ent}.
If $h\in C_\nu(X)$, the
function $A_a(t)h$ is supported in $T^a\{ y\in B\tq \phi(y)>a\}$,
whose measure is $m[\phi>a] =O(1/a^\beta)$. Thus, $\left| \int_X
A_a(t) h \right| \leq \frac{C}{a^\beta} \norm{h}_\infty \leq
\frac{C}{a^{\beta-1}} \norm{h}$. Moreover, $\norm{E(k,t)
B_b(t) (u)} \leq e(k,t) \norm{B_b(t)} \leq e(k,t)
\frac{C}{b^{\beta-1}}$ by Lemma \ref{int_Bb}. Thus,
  \begin{equation}
  \label{majore_deuxieme_somme}
  \begin{split}
  \left|\sum_{a+k+b=n} \int_X A_a(t) E(k,t) B_b(t) (u)v\right|&
  \leq \sum_{a+k+b=n} \frac{C}{a^{\beta-1}} e(k,t)
  \frac{C}{b^{\beta-1}}
  \\&
  \leq C \left(\frac{1}{n^{\beta-1}} \right) \star e(n,t) \star
  \left( \frac{1}{n^{\beta-1}} \right)
  \leq e(n,t)
  \end{split}
  \end{equation}
by \eqref{convoles_ent}.

The first sum of \eqref{decompose_somme_avec_Ent} can be written
as
  \begin{equation*}
  I=
  \sum_{a+k+b=n} \left(\int_X A_a(t)(1_B)v \right)
  \left(1- \frac{\sigma^2}{2}L(t)\right)^k \left( \int_B B_b(t)(u)
  \right).
  \end{equation*}
Using Lemmas \ref{int_Aa} and \ref{int_Bb} and convolving, we find
a sequence $w_n$ such that
$w_n=O(1/n^\beta)$, $\sum w_n=\left(\int u\right)\left(\int v\right)$, and
  \begin{equation*}
  I= \left( w_n +O\left(\frac{|t|}{n^{\beta-1}} \right) \right)
  \star \left(1- \frac{\sigma^2}{2}L(t)\right)^k.
  \end{equation*}
As $L(t)\sim t^2$ when $t\to 0$, the term coming from
$O\left(\frac{|t|}{n^{\beta-1}} \right)$ is bounded by
$|t| \left(\frac{1}{n^{\beta-1}}\right)\star (1-dt^2)^n \leq
e(n,t)$. Moreover, for $x,y\in \R$,
  \begin{equation}
  \label{estime_exp}
  |e^x-e^y| \leq |x-y| e^{\max(x,y)}.
  \end{equation}
%\begin{proof}
%On peut supposer que $a\geq b$. On factorise alors par $e^b$, puis
%on applique l'in\'egalit\'e des accroissements finis, qui donne
%$|e^x-1|\leq xe^x$ quand $x\geq 0$.
%\end{proof}
Thus,
  \begin{align*}
  \left|
  \sum_{k=0}^n w_{n-k} \left(1- \frac{\sigma^2}{2}L(t)\right)^k
  -\sum_{k=0}^n w_{n-k} \left(1- \frac{\sigma^2}{2}L(t)\right)^n
  \right|
  \!\!\!\!\!\!\!\!\!\!\!\!\!\!\!\!\!\!\!
  \!\!\!\!\!\!\!\!\!\!\!\!\!\!\!\!\!\!\!
  \!\!\!\!\!\!\!\!\!\!\!\!\!\!\!\!\!\!\!
  %\!\!\!\!\!\!\!\!\!\!\!\!\!\!\!\!\!\!\!
  %\!\!\!\!\!\!\!\!\!\!\!\!\!\!\!\!\!\!\!
  %\!\!\!\!\!\!\!\!\!\!\!\!\!\!\!\!\!\!\!
  &
  \\&
  \leq \sum_{k=0}^n \frac{C}{(n-k)^\beta} (n-k)
  \left| \ln\left(1- \frac{\sigma^2}{2}L(t)\right) \right|
  \left(1- \frac{\sigma^2}{2}L(t)\right)^k
  \\&
  \leq \sum_{k=0}^n C t^2 (1-d t^2)^{k}
  \frac{1}{(n-k)^{\beta-1}}
  %\leq \frac{C}{n^{\beta-1}}+C t^2
  %\left(1-\frac{d}{2}t^2\right)^n
  \leq  e(n,t).
  \end{align*}

Hence, up to $O(e(n,t))$, the integral $\int e^{it S_n f} \cdot u
\cdot v\circ T^n$ is
equal to $\sum_{k=0}^n w_{n-k} \left(1-\frac{\sigma^2}{2}L(t)
\right)^n$. Finally,
   \begin{multline*}
   \left|\left(1-\frac{\sigma^2}{2}L(t) \right)^n \int_X u \int_X v-
   \sum_{k=0}^n w_{n-k}
   \left(1-\frac{\sigma^2}{2}L(t) \right)^n \right|\\
   =\left(1-\frac{\sigma^2}{2}L(t) \right)^n \left|\sum_{k=n+1}^\infty
   w_k\right|
   \leq C\sum_{k=n+1}^\infty \frac{1}{k^\beta}
   \leq \frac{C}{n^{\beta-1}} \leq e(n,t).
   \end{multline*}
This concludes the proof.
\end{proof}

When $u=v=1$, Theorem \ref{estimee_clef_markov} states that, for small
$t$, the
characteristic function of $S_n f$ behaves essentially like
$\left(1-\frac{\sigma^2}{2}L(t)\right)^n$, which is very similar
to the characteristic
function of the sum of $n$ independent identically distributed random
variables. Hence, using this estimate, it will be possible to use the classical
probabilistic proofs to get the local limit theorem or the
Berry-Esseen theorem. However, some care is still required to check
that the error term in Theorem \ref{estimee_clef_markov} is
sufficiently small so that these proofs still work.

\section{Proof of the local limit theorem}

\label{section_limite_locale}

\subsection{Periodicity problems}

\label{subsection_periodicite}
This paragraph is related to the end of Section 3 of
\cite{aaronson_denker}. The differences come mainly from the inductive
process and the fact that we are considering series of operators
instead of a single operator.

Let $X$ be a Young tower with $\gcd(\phi_i)=1$, $\tau<1$ and $f\in
C_\tau(X)$. For $t\in \R$, $e^{itf}$ is said to be
\emph{cohomologous} to a constant $\lambda \in S^1$ if there exists
$\omega:X\to S^1$ such that $e^{itf} =\lambda \frac{\omega\circ
T}{\omega}$ almost everywhere.

\begin{prop}
\label{prop_aperiodicite}
Let $t\in \R$, and $z\in S^1$. The following assertions are equivalent:
\begin{enumerate}
\item $e^{itf}$ is cohomologous to $z^{-1}$.
\item $1$ is an eigenvalue of the operator $R(z,t)$ acting on
$C_\nu(B)$.
\end{enumerate}
\end{prop}
\begin{proof}
Suppose first that there exists a nonzero $\omega\in C_\nu(B)$
such that $R(z,t)\omega=\omega$. As $R(z,t)(v)=\hat{T}_B (z^\phi
e^{it f_B} v)$, the adjoint operator of $R(z,t)$ in $L^2$ is
$W(v)=z^{-\phi} e^{-it f_B} v\circ T_B$. Then
  %\begin{align*}
  $\norm{W\omega-\omega}_{L^2}^2
  %&= \norm{Wg}_{L^2}^2 - 2\real\langle Wg,g\rangle + \norm{g}_{L^2}^2
  %=\norm{Wg}_{L^2}^2 -2\real \langle g, R(z,t)g \rangle +\norm{g}_{L^2}^2
  %\\&
  %=\norm{Wg}_{L^2}^2 -2\real \langle g,g \rangle +\norm{g}_{L^2}^2
  =\norm{W\omega}_{L^2}^2 - \norm{\omega}_{L^2}^2.$
  %\end{align*}
As $T_B$ preserves the measure, we have
$\norm{\omega}_{L^2}=\norm{W\omega}_{L^2}$, whence
$\norm{W\omega-\omega}_{L^2}=0$. Consequently, $\omega=z^{-\phi}
e^{-it f_B} \omega\circ T_B$ almost everywhere. Taking the
modulus, $|\omega|=|\omega|\circ T_B$, and the ergodicity of $T_B$
gives that $|\omega|$ is almost everywhere constant. We can assume
that $|\omega|=1$. We extend the function $\omega$ to $X$ by
setting
  \begin{equation*}
  \omega(x,k) = \omega(x,0) e^{itf(x,0)} \cdots e^{itf(x,k-1)} z^k.
  \end{equation*}
Then, for $k<\phi(x)-1$, we have by construction $\omega\circ
T(x,k)/\omega(x,k)= ze^{itf(x,k)}$. Moreover, for $k=\phi(x)-1$,
  \begin{equation*}
  \omega \circ T(x,k)/\omega(x,k)
  =\omega\circ T_B(x,0)/\omega(x,k)
  =e^{it f_B(x,0)}z^{k+1} \omega(x,0)/ \omega(x,k)
  =e^{it f(x,k)} z.
  \end{equation*}
Thus, $e^{itf} = z^{-1} \omega \circ T/ \omega$ almost everywhere.

Conversely, suppose that a measurable function $\omega$ satisfies $e^{itf}
= z^{-1} \omega\circ T/\omega$. The previous calculations give
$e^{it f_B} =z^{-\phi}\omega\circ T_B /\omega$. The operator
$R(z,t)=\hat{T}_B(z^{\phi}e^{itf_B} \cdot)$ acts on $L^1$, and
satisfies
  \begin{equation*}
  R(z,t) (\omega)=\hat{T}_B\left( z^\phi e^{it f_B} \omega\right)
  =\hat{T}_B \left(\omega\circ T_B\right)
  =\omega.
  \end{equation*}
But, in Lemma
\ref{enumere_estimees}, we have seen that $R(z,t)$ satisfies a
Doeblin-Fortet inequality between the spaces $C_\nu(B)$ and
$L^1(B)$. \cite{ionescu_tulcea_marinescu} ensures that the
eigenfunctions of $R(z,t)$ in $L^1(B)$ for the eigenvalue $1$  are
in fact in $C_\nu(B)$, i.e.\ $\omega\in C_\nu(B)$.
\end{proof}

\begin{cor}
\label{aperiodicite}
The set $\mathfrak{A}:=\{t\in \R \tq e^{itf} \text{ is cohomologous to
a constant}\}$ is a closed subgroup of $\R$. Moreover, for every $t\in
\mathfrak{A}$, there exists a unique $z(t)\in S^1$ such that $e^{itf}$
is cohomologous to $z(t)$. Finally, the map $t \mapsto z(t)$ is a
continuous morphism from $\mathfrak{A}$ to $S^1$.
\end{cor}
\begin{proof}
The set $\mathfrak{A}$ is clearly a subgroup of $\R$. If $e^{itf}$ is
cohomologous simultaneously to $z$ and $z'$, then $z'=\frac{e^{ih
\circ T}}{e^{ih}} z$ for some function $h:X\to \R$. As $T$ is mixing
(\cite[Theorem 1 (iii)]{lsyoung:recurrence}), the only constant $s$
satisfying $g\circ T=sg$ for some nonzero function $g$ is $1$. This
implies that $z=z'$. The map $t\mapsto z(t)$ is thus well defined,
and it is clearly a group morphism.

It remains to check that $\mathfrak{A}$ is closed and that $z(t)$
is continuous. Let $t_n$ be a sequence of $\mathfrak{A}$
converging to $T \in \R$. Let $Z$ be a cluster point of the
sequence $z(t_n)$. By \ref{point_5_lemme} in Lemma
\ref{enumere_estimees}, $(z,t)\mapsto R(z,t)$ is a continuous map
with values in $\Hom(C_\nu(B))$. If $I-R(Z^{-1},T)$ were
invertible, then $I-R(z_n^{-1},t_n)$ would also be invertible for
large enough $n$, which is a contradiction by the previous
proposition. Thus, $I-R(Z^{-1},T)$ is not invertible, whence $1$
is an eigenvalue of $R(Z^{-1},T)$ by quasi-compactness. This
implies that  $T \in \mathfrak{A}$ and $Z=z(T)$, once again by the
previous proposition.
\end{proof}

Consequently, there are three cases to be considered for the local
limit theorem: $\mathfrak{A}$ is either $\R$, or $\{0\}$, or a
discrete subgroup of $\R$. If it is $\R$, \cite{moore_schmidt}
ensures that $f$ can be written as $g-g\circ T$, hence $\sigma=0$ in
the central limit theorem, and there is nothing to prove. If
$\mathfrak{A}=\{0\}$, it is not possible to write $f$ as
$\rho+g-g\circ T+\lambda q$, where $\rho\in \R$, $g:X\to \R$ is
measurable, $\lambda >0$ and $q:X \to \Z$, i.e.\ $f$ is aperiodic.
This case is handled by Theorem \ref{limite_locale_abstrait}.
Finally, $\mathfrak{A}=2\pi\Z$ means that $f=\rho+g-g\circ T+q$ where
$q$ takes integer values, and that there is no such expression
where $q$ takes its values in $n\Z$ with $n\geq 2$. This is dealt
with in Theorem \ref{limite_locale_periodique_abstrait}.

\subsection{The aperiodic case}

\begin{proof}[Proof of Theorem
\ref{limite_locale_abstrait}]

Let $f$ be an aperiodic function on $X$, as in the hypotheses of
Theorem \ref{limite_locale_abstrait}. Then, for every $(z,t)\in
(\overline{\D}\times \R) -\{(1,0)\}$, $I-R(z,t)$ is invertible: for
$|z|<1$, the spectral radius of $R(z,t)$ is at most $|z|<1$
(\ref{point_2_lemme} in Lemma \ref{enumere_estimees}), and for
$|z|=1$ this comes from Proposition \ref{prop_aperiodicite} and
Corollary \ref{aperiodicite}, since $\mathfrak{A}=\{0\}$.

Let $\alpha>0$ be given by Theorem \ref{estimee_clef_markov}: we
control the behavior of the integrals when $|t|\leq \alpha$. Take
$K>0$. We will show the following fact: \emph{there exists $C>0$
such that, for every $|t|\in [\alpha,K]$, for every $n\in \N^*$,
for all functions $u\in C_\tau(X)$ and $v\in L^\infty(X)$, }
  \begin{equation}
  \label{int_sur_t_grand}
  \left|\int_X e^{it S_n f}\cdot u \cdot v \circ T^n\right|
  \leq \frac{C}{n^{\beta-1}}\norm{u} \norm{v}_\infty.
  \end{equation}
Let us write $\boA=\boO_{\beta-1}^+(\Hom(C_{\nu}(B)))$ (the Banach
algebra of series $\sum_{0}^\infty A_n z^n$ where $A_n \in
\Hom(C_{\nu}(B))$ and $\norm{A_n}=O(1/n^{\beta-1})$, with the norm
$\norm{\sum A_n z^n} =\sup_{n\in \N} (n+1)^{\beta-1} \norm{A_n}$). The map
$t\mapsto R(z,t)$ is continuous from $[-K,-\alpha]\cup [\alpha,K]$
to $\boA$ (\ref{point_5_lemme} in Lemma \ref{enumere_estimees}).
Moreover, $I-R(z,t)$ is invertible on $\overline{\D}$ for $t$ in
these intervals. Theorem \ref{wiener_D} shows that
$(I-R(z,t))^{-1} \in \boA$, and the continuity of the inversion
even yields that $t\mapsto (I-R(z,t))^{-1}$ is continuous. By
compactness, there exists $C$ such that
$\norm{(I-R(z,t))^{-1}}_{\boA} \leq C$ for $|t|\in [\alpha,K]$. As
$(I-R(z,t))^{-1} =\sum T_n(t) z^n$, this implies that
$\norm{T_n(t)} \leq \frac{C}{n^{\beta-1}}$, uniformly in $n$ and
$t$.

We have
  \begin{equation*}
  \int_X e^{it S_n f}\cdot u \cdot v\circ T^n
  =\int_X C_n(t)(u)v + \sum_{a+k+b=n} \int_X A_a(t) T_k(t) B_b(t)(u)v.
  \end{equation*}
By Lemma \ref{int_Cn}, $ \left|\int_X C_n(t)(u)v\right| \leq
\frac{C}{n^{\beta-1}} \norm{u}_\infty \norm{v}_\infty$.
By Lemmas \ref{int_Aa} and \ref{int_Bb},
  \begin{align*}
  \left|\int A_a(t) T_k(t) B_b(t)(u)v \right| &\leq
  \frac{C}{a^{\beta-1}} \norm{T_k(t) B_b(t)(u)}_\infty \norm{v}_\infty
  \leq \frac{C}{a^{\beta-1}} \norm{T_k(t)} \norm{B_b(t)} \norm{u}
  \norm{v}_\infty
  \\&
  \leq \frac{C}{a^{\beta-1}} \frac{C}{k^{\beta-1}}
  \frac{C}{b^{\beta-1}} \norm{u}\norm{v}_\infty.
  \end{align*}
Thus,
  \begin{equation*}
  \left|\int_X e^{it S_n f}\cdot u\cdot v\circ T^n \right|
  \leq C \left(\frac{1}{n^{\beta-1}} \right) \star
  \left(\frac{1}{n^{\beta-1}} \right) \star
  \left(\frac{1}{n^{\beta-1}} \right) \norm{u}\norm{v}_\infty
  \leq \frac{C}{n^{\beta-1}} \norm{u}\norm{v}_\infty,
  \end{equation*}
and \eqref{int_sur_t_grand} is proved.

We prove now the local limit theorem, using the method of Breiman
(\cite{breiman}). Take $u,v$ and $k_n$ as in the assumptions of Theorem
\ref{limite_locale_abstrait}. Let $\psi \in L^1(\R)$ be such that its
Fourier transform $\hat{\psi}$ is supported in
$[-K,K]$. Then $\psi(x)=\frac{1}{2\pi} \int_{-K}^K
\hat{\psi}(t)e^{itx} \dd t$, whence
  \begin{align}
  \label{equation_preuve_locale}
  \raisetag{58pt}
  \begin{split}
  \sqrt{n} E(\psi(S_n f-k_n-u-v\circ T^n))
  &
  =\frac{\sqrt{n}}{2\pi} \int_{-K}^K
  \hat{\psi}(t) E\left(e^{it (S_n f-k_n-u-v\circ T^n)}\right) \dd t
  \\&
  =\frac{\sqrt{n}}{2\pi} \int_{-\alpha}^\alpha \hat{\psi}(t) e^{-it
  k_n} E(e^{it S_n f}e^{-itu}e^{-it v\circ T^n}) \dd t
  \\&\hphantom{=\ }
  +\frac{\sqrt{n}}{2\pi} \int_{\alpha \leq |t| \leq K}
  \hat{\psi}(t) e^{-it
  k_n} E(e^{it S_n f}e^{-itu}e^{-it v\circ T^n}) \dd t.
  \end{split}
  \end{align}
For $\alpha \leq |t| \leq K$, the norms $\norm{ e^{-itu}}$ and
$\norm{e^{-itv}}_\infty$ remain bounded. Hence,
\eqref{int_sur_t_grand} implies that
$\left|E(e^{it S_n f}e^{-itu}e^{-it v\circ T^n})
\right|\leq \frac{C}{n^{\beta-1}}$. Therefore, the
second integral tends to $0$. For
the first integral, we approximate $E(e^{it S_n f}e^{-itu}e^{-itv\circ
T^n})$ by
$\left(1-\frac{\sigma^2}{2}L(t) \right)^n \int e^{-itu}\int e^{-itv}$.
By Theorem \ref{estimee_clef_markov}, the error term is bounded by
  \begin{equation*}
  C\sqrt{n}
  \int_{-\alpha}^\alpha\left( \frac{1}{n^{\beta-1}} +|t|
  \left(\frac{1}{n^{\beta-1}} \right) \star (1-dt^2)^n\right)\dd t.
 \end{equation*}
Let us
show that this integral tends to $0$. This is clear for the first
term. For the second term, we cut the integral in two pieces. For
$|t|\leq \frac{1}{\sqrt{n}}$, the convolution is bounded (since
$\frac{1}{n^{\beta-1}}$ is summable and $(1-dt^2)^n \leq 1$), whence
the integral is $\leq C \sqrt{n}  \int_{|t|\leq 1/\sqrt{n}} |t|\dd t \to
0$. For $|t|\geq \frac{1}{\sqrt{n}}$, Lemma \ref{ngamma_star_t} gives
that the convolution is bounded by $\frac{1}{|t|n^{\beta-1}}+|t|
(1-dt^2)^n$, whence the integral is less than
  \begin{equation*}
  C \sqrt{n} \int_{1/\sqrt{n} \leq |t| \leq \alpha }
  \frac{1}{|t|n^{\beta-1}}+|t| (1-dt^2)^n\dd t
  \leq C\sqrt{n} \frac{\ln n}{n^{\beta-1}}
  +C \sqrt{n} \left[ \frac{(1-dt^2)^{n+1}}{-2d(n+1)}
  \right]_{1/\sqrt{n}}^\alpha
  \to 0.
  \end{equation*}
Finally, we have proved that
  \begin{multline*}
  \sqrt{n} E(\psi(S_n f-k_n-u-v\circ T^n))
  \\
  =\frac{\sqrt{n}}{2\pi} \int_{-\alpha}^\alpha \hat{\psi}(t)
  e^{-it k_n} \left(1-\frac{\sigma^2}{2}
  L(t) \right)^n E(e^{-itu})E(e^{-itv})\dd t  +o(1).
  \end{multline*}
But
  \begin{align*}
  \frac{\sqrt{n}}{2\pi} \int_{-\alpha}^\alpha \hat{\psi}(t)
  e^{-it k_n}\left(1-\frac{\sigma^2}{2}
  L(t) \right)^n E(e^{-itu})E(e^{-itv}) \dd t
  \!\!\!\!\!\! \!\!\!\!\!\! \!\!\!\!\!\! \!\!\!\!\!\!
  \!\!\!\!\!\! \!\!\!\!\!\! \!\!\!\!\!\! \!\!\!\!\!\!
  \!\!\!\!\!\! \!\!\!\!\!\! \!\!\!\!\!\! \!\!\!\!\!\!
  \!\!\!\!\!\! \!\!\!\!\!\! \!\!\!\!\!\! \!\!\!\!\!\!
   \!\!\!\!\!\!
  \\&
  =\frac{1}{2\pi} \int_{-\alpha \sqrt{n}}^{\alpha \sqrt{n}}
  \hat{\psi}\left(\frac{t}{\sqrt{n}} \right) e^{-it k_n/\sqrt{n}}
  \left(1-\frac{\sigma^2}{2} L\left(\frac{t}{\sqrt{n}}\right)
  \right)^n  E(e^{-i \frac{t}{\sqrt{n}} u})E(e^{-i\frac{t}{\sqrt{n}} v}) \dd t
  \\&
  \to \frac{1}{2\pi} \int_\R \hat{\psi}(0) e^{-it \kappa}
  e^{-\frac{\sigma^2}{2}t^2} \dd t,
  \end{align*}
by dominated convergence. We have used the fact that $L(t)\sim t^2$
close to $0$, and in particular, if $\alpha$ is small enough,
$\left(1-\frac{\sigma^2}{2} L \left(\frac{t}{\sqrt{n}}\right)
\right)^n \leq \left(1-\frac{\sigma^2 t^2}{4 n} \right)^n \leq
e^{-\frac{\sigma^2 t^2}{4}}$, which gives the domination.

Set $\chi(\kappa)=\frac{e^{-\frac{\kappa^2}{2\sigma^2}}}{\sigma
\sqrt{2\pi}}$. We have proved that, for any $\psi$ in $L^1$ with
$\hat{\psi}$ compactly supported,
  \begin{equation}
  \label{But_thm_limite_locale}
  \sqrt{n} E( \psi(S_n f-k_n-u-v\circ T^n))
  \to \chi(\kappa) \int_\R \psi(x)\dd x.
  \end{equation}
Equation \eqref{But_thm_limite_locale} can then be extended to a
larger class of functions by density arguments (see \cite{breiman}),
and this larger class contains in particular the characteristic
functions of bounded intervals. This concludes the proof.
\end{proof}

\subsection{The periodic case}

The following theorem gives the local limit theorem when the group
$\mathfrak{A}$ of Paragraph \ref{subsection_periodicite} is a
discrete subgroup of $\R$, for example $2\pi\Z$.

\begin{thm}[local limit theorem, periodic case]
\label{limite_locale_periodique_abstrait}
Let $X$ be a Young tower with $\gcd(\phi_i)=1$. Assume that
$m[\phi>n]=O(1/n^{\beta})$ with $\beta>2$. Let $\tau<1$.
Take $f\in C_\tau(X)$ of
zero integral, and $\sigma^2$ given by Theorem \ref{TCL_normal}.

Assume that $f=\rho+q$ where $q$ takes integer values and $\rho \in
\R$, but that $f$ can not be written as $f=\rho'+g-g\circ T+\lambda
q'$, where $\lambda \in \N-\{1\}$ and $q':X \to \Z$ (this implies
in particular $\sigma>0$). Then, for every sequence $k_n$ with
$k_n-n\rho \in \Z$ such that $k_n / \sqrt{n} \to \kappa\in \R$,
  \begin{equation*}
  \sqrt{n}\; m\left\{ x \in X \tq S_n f(x)=k_n\right\}
  \to \frac{e^{-\frac{\kappa^2}{2\sigma^2}}}{\sigma \sqrt{2\pi}}.
  \end{equation*}
\end{thm}
\begin{proof}
This is essentially the same proof as that of Theorem
\ref{limite_locale_abstrait}, but we use the Fourier transform on
$\Z$ (i.e.\ Fourier series) instead of the Fourier transform on
$\R$.

If $k$ and $l$ are two integer numbers,
  \begin{equation*}
  1_{k=l} = \frac{1}{2\pi} \int_{-\pi}^{\pi} e^{it (l-k)} \dd t.
  \end{equation*}
Applying this equation to $k_n-n\rho$ and $S_n f(x)-n \rho$ and
integrating gives
  \begin{equation*}
  m\left\{ x\in X \tq S_n f(x)=k_n \right\}
  =
  \frac{1}{2\pi} \int_{-\pi}^{\pi} e^{-itk_n} E( e^{it S_n f}) \dd t.
  \end{equation*}
This is an analogue of \eqref{equation_preuve_locale}. From this
point on, the proof of Theorem \ref{limite_locale_abstrait}
applies. The only problem is to check that, on $[-\pi,-\alpha]\cup
[\alpha,\pi]$, $I-R(z,t)$ is invertible for $z\in\overline{\D}$.
This comes from the assumptions on $f$, which ensures that
$I-R(z,t)$ is invertible as soon as $t\not\in 2\pi \Z$, by
Proposition \ref{prop_aperiodicite} and Corollary
\ref{aperiodicite}, since $\mathfrak{A}=2\pi\Z$.
\end{proof}

\section{Proof of the central limit theorem with speed}

\begin{proof}[Proof of Theorem \ref{TCL_vitesse}]

The Berry-Esseen Theorem (\cite{feller:2}) implies that the result will be
proved if we show that, for some $c>0$,
  \begin{equation*}
  \int_{-c\sqrt{n}}^{c\sqrt{n}}\frac{1}{|t|} \left|
  E(e^{i\frac{t}{\sqrt{n}}S_n f}) -e^{-\frac{\sigma^2}{2}t^2}
  \right|\dd t =O\left(\frac{1}{n^{\delta/2}}\right).
  \end{equation*}

We first estimate the integral between $-1/n$ and $1/n$:
  \begin{align*}
  \int_{-1/n}^{1/n} \frac{1}{|t|}\left|E(e^{i\frac{t}{\sqrt{n}}S_n f})
  -e^{-\frac{\sigma^2}{2}t^2}\right| \dd t
  &\leq \int_{-1/n}^{1/n} \frac{1}{|t|} \left|E(e^{i\frac{t}{\sqrt{n}}S_n f})
  -1\right|\dd t
  +\int_{-1/n}^{1/n} \frac{1}{|t|} \left|e^{-\frac{\sigma^2}{2} t^2}-1
  \right|\dd t
  \\&
  \leq \int_{-1/n}^{1/n} \frac{1}{\sqrt{n}}E(|S_n f|)\dd t
  +\int_{-1/n}^{1/n} \frac{\sigma^2}{2}|t|\dd t.
  \end{align*}
But $\int |S_n f|\leq n\int |f|$, whence we get $O(1/\sqrt{n})$ for
this term.

Let $L(t)$ be given by Proposition \ref{TCL_induit}. Then,
for small enough $c$,
  \begin{align*}
  \int_{1/n\leq |t| \leq c\sqrt{n}} \frac{1}{|t|}
  \left|E(e^{i\frac{t}{\sqrt{n}}S_n f})
  -e^{-\frac{\sigma^2}{2}t^2}\right|\dd t
  \!\!\!\!\!\!\!\!\!\!\!\!\!\!\!
  \!\!\!\!\!\!\!\!\!\!\!\!\!\!\!\!\!\!\!\!
  &
  \\&
  \leq
  \int_{1/n\leq |t| \leq c\sqrt{n}} \frac{1}{|t|}\left|
  \left(1-\frac{\sigma^2}{2}
  L\left(\frac{t}{\sqrt{n}}\right)\right)^n
  -e^{-\frac{\sigma^2}{2}t^2}\right|\dd t
  \\&\ \ +
  \int_{1/n\leq |t| \leq c\sqrt{n}} \frac{1}{|t|}
  \left|E(e^{i\frac{t}{\sqrt{n}}S_n f})
  -\left(1-\frac{\sigma^2}{2}L\left(\frac{t}{\sqrt{n}}\right)\right)^n
  \right|\dd t.
  \end{align*}
Let us show that the second term satisfies
  \begin{equation*}
  \int_{1/n\leq |t| \leq c\sqrt{n}} \frac{1}{|t|}
  \left|E(e^{i\frac{t}{\sqrt{n}}S_n f})
  -\left(1-\frac{\sigma^2}{2} L\left(
  \frac{t}{\sqrt{n}}\right)\right)^n \right|\dd t =
  O\left(\frac{1}{\sqrt{n}}\right).
  \end{equation*}

Set $e(n,t)=C\left[\frac{1}{n^{\beta-1}}
+|t|\left(\frac{1}{n^{\beta-1}}\right) \star (1-dt^2)^n \right]$.
By Theorem \ref{estimee_clef_markov},
  \begin{equation*}
  \left|E(e^{i\frac{t}{\sqrt{n}}S_n f}) -\left(1-\frac{\sigma^2}{2}
  L\left( \frac{t}{\sqrt{n}}\right)\right)^n\right| \leq e(n,
  t /\sqrt{n}).
  \end{equation*}
Thus, it is enough to prove that
  \begin{equation*}
  \int_{1/n \leq |t| \leq c\sqrt{n}} \frac{e(n,t/\sqrt{n})}{|t|}\dd t =
  O\left(\frac{1}{\sqrt{n}} \right).
  \end{equation*}

For the term $\frac{1}{n^{\beta-1}}$in $e(n,t)$, the integral
is $C\frac{\ln n}{n^{\beta-1}}$, which is $O(1/\sqrt{n})$.

For the other term $c(n,t)=|t| \left(\frac{1}{n^{\beta-1}}\right)
\star (1-dt^2)^n$, we cut the integral in two pieces. For $|t|\leq
1$, the convolution is bounded (since $\frac{1}{n^{\beta-1}}$ is
summable and $(1-d\frac{t^2}{n})^n \leq 1$). It remains $\int_{1/n}^1
\frac{1}{|t|} \frac{|t|}{\sqrt{n}}\dd t\leq \frac{1}{\sqrt{n}}$. For $|t|
\geq 1$, we use Lemma \ref{ngamma_star_t}, which gives that
 $c(n,t)\leq
\frac{1}{|t|n^{\beta-1}}+ |t|\left(1-\frac{d}{2}t^2\right)^n$.
Thus,
  \begin{equation*}
  \int_1^{\sqrt{n}} \frac{1}{|t|}
  \left|c\left(n,\frac{t}{\sqrt{n}}\right)\right|\dd t
  \leq \int_1^{\sqrt{n}} \frac{\sqrt{n}}{t^2 n^{\beta-1}}\dd t
  + \frac{1}{\sqrt{n}}
  \int_1^{\sqrt{n}} e^{-dt^2/2} \dd t
  =O\left(\frac{1}{\sqrt{n}} \right).
  \end{equation*}

Finally, we have proved that
  \begin{multline}
  \label{uziyeroyaiozrye}
  \int_{-c\sqrt{n}}^{c\sqrt{n}}\frac{1}{|t|} \left|
  E(e^{i\frac{t}{\sqrt{n}}S_n f}) -e^{-\frac{\sigma^2}{2}t^2}
  \right|\dd t \\
  \leq \int_{|t| \leq c\sqrt{n}} \frac{1}{|t|}\left|
  \left(1-\frac{\sigma^2}{2}
  L\left(\frac{t}{\sqrt{n}}\right)\right)^n
  -e^{-\frac{\sigma^2}{2}t^2}\right|\dd t +
  O\left(\frac{1}{\sqrt{n}}\right).
  \end{multline}

We have only to deal with the powers of a function. Hence, it will be
possible to use the same methods as in probability theory. More
precisely, the study of the speed in the central limit theorem in
\cite[Theorem
3.4.1]{ibragimov_linnik} uses the two following facts:
\begin{enumerate}
\item
\label{lkjsqdflkmjq}
If a random variable $Z$ satisfies $E(|Z|^2
1_{|Z|>z})= O(z^{-\delta})$ with $0<\delta\leq 1$ and,
in the $\delta=1$ case, $E(Z^3 1_{|Z|\leq z})=O(1)$, then there exists
a constant
$\lambda^2\geq 0$ such that
$E(e^{itZ})=1-\frac{\lambda^2}{2} t^2(1+\gamma(t))$, with
$\int_{-x}^x t^2 |\gamma(t)|\dd t=O(x^{3+\delta})$ when $x \to 0$.

\item
\label{jklqhsdfjklhuiyi}
If a function $\tilde{\gamma}(t)$ satisfies $\int_{-x}^x t^2
|\tilde{\gamma}(t)|\dd t = O(x^{3+\delta})$ with $0<\delta \leq 1$, and
$\sigma^2 \geq 0$, then $\Lambda(t):= t^2(1+\tilde{\gamma}(t))$
satisfies $\int_{|t| \leq c\sqrt{n}}
\frac{1}{|t|}\left|
  \left(1-\frac{\sigma^2}{2} \Lambda \left(\frac{t}{\sqrt{n}}\right)
  \right)^n
  -e^{-\frac{\sigma^2}{2}t^2}\right|\dd t = O(n^{-\delta/2})$.
\end{enumerate}

By Proposition \ref{TCL_induit}, the eigenvalue $\lambda(t)$ of
$\hat{T}_B(t)$  is equal to $E_B(e^{itf_B})+\alpha t^2 +O(t^3)$ for
some constant $\alpha$. The fact  \ref{lkjsqdflkmjq} applied to the
random variable $f_B:B\to \R$ implies that
$E_B(e^{itf_B})=1-\frac{\lambda^2}{2}t^2(1+\gamma(t))$ for some
function $\gamma(t)$ satisfying $\int_{-x}^x t^2
|\gamma(t)|\dd t=O(x^{3+\delta})$. As
$\lambda(t)=1-\frac{\sigma^2}{2m(B)}L(t)$, this implies that
$L(t)=t^2\bigl(1+m(B)\frac{\lambda^2}{\sigma^2}\gamma(t)
+O(t)\bigr)$. Hence, we can write
$L(t)=t^2(1+\tilde{\gamma}(t))$ with $\int_{-x}^x t^2
|\tilde{\gamma}(t)|\dd t=O(x^{3+\delta})$. Therefore, the fact
\ref{jklqhsdfjklhuiyi} implies that $\int_{|t| \leq c\sqrt{n}}
\frac{1}{|t|}\left| \left(1-\frac{\sigma^2}{2}
L\left(\frac{t}{\sqrt{n}}\right)\right)^n
-e^{-\frac{\sigma^2}{2}t^2}\right|\dd t=O(n^{-\delta/2})$.

By \eqref{uziyeroyaiozrye}, we obtain
$\int_{-c\sqrt{n}}^{c\sqrt{n}}\frac{1}{|t|} \left|
  E(e^{i\frac{t}{\sqrt{n}}S_n f}) -e^{-\frac{\sigma^2}{2}t^2}
  \right|\dd t=O(n^{-\delta/2})$, which concludes the proof.
\end{proof}

\appendix

\section{The Wiener Lemma}
\label{appendice_preuve_wiener}

In this appendix, we prove that the algebra $\boO_\gamma(\boC)$
introduced in Paragraph \ref{subsection_wiener} is indeed a Banach
algebra, and that it satisfies a Wiener Lemma (Theorem
\ref{wiener_S1}).

Let $\boC$ be a Banach algebra, and take $\gamma>1$.
Write  $w_n=(n+1)^{-\gamma}$ for $n\geq 0$. There exists a
constant $c$ such that $(w_n) \star (w_n) \leq c w_n$. We define a
norm on $\boO_\gamma(\boC)$ by
  \begin{equation}
  \label{norme_Ogamma}
  \norm{\sum_{n\in \Z} A_n z^n}
  =\left(\sum_{n\in \Z} \norm{A_n} + c \sup_{n\geq 0}
  \frac{\norm{A_n}}{w_n} \right)
  +
  \left(\sum_{n\in \Z} \norm{A_n} + c \sup_{n\leq 0}
  \frac{\norm{A_{n}}}{w_{|n|}} \right).
  \end{equation}

\begin{prop}
\label{wiener_algebra}
Let $\boC$ be a Banach algebra, and $\gamma>1$.
With the norm \eqref{norme_Ogamma}, $\boO_\gamma(\boC)$ is a Banach algebra.
\end{prop}
\begin{proof}
The completeness is clear. It is sufficient to prove the
submultiplicativity of the norm for one half of this norm, for
example the first one. Let us write $\norm{ \sum A_n z^n}_1= \sum
\norm{A_n}$, and $P_w(\sum A_n z^n) =\sup_{n\geq 0}
\frac{\norm{A_n}}{w_n}$. Then, if  $A=\sum A_n z^n$ and $B=\sum B_n
z^n$, we have $\norm{AB}_1 \leq \norm{A}_1 \norm{B}_1$. Moreover, for
$n\geq 0$,
  \begin{align*}
  \frac{\norm{(AB)_n}}{w_n}&
  \leq \frac{\sum_k \norm{A_k B_{n-k}}}{w_n}
  \leq \frac{1}{w_n}\sum_{k=0}^n \norm{A_k B_{n-k}}
  +\left(\sum_{k=-\infty}^{-1}\norm{A_k} P_w(B)
\frac{w_{n-k}}{w_n}\right)
  \\&
  \hphantom{\leq \frac{\sum_k \norm{A_k B_{n-k}}}{w_n}\leq \ }
  +\left(\sum_{k=-\infty}^{-1} \norm{B_k} P_w(A) \frac{w_{n-k}}{w_n}
  \right)
  \\&
  \leq P_w(A)P_w(B) \frac{(w\star w)_n}{w_n} +\norm{A}_1 P_w(B)
  +\norm{B}_1 P_w(A).
  \end{align*}
Thus, $P_w(AB)\leq cP_w(A)P_w(B)+\norm{A}_1 P_w(B)+ \norm{B}_1
P_w(A)$. This gives the conclusion.
\end{proof}

We will now identify the characters of the commutative algebra
$\boO_\gamma(\C)$, i.e.\ the algebra morphisms from
$\boO_\gamma(\C)$ to $\C$. For $\lambda \in S^1$, we can define a
character $\chi_\lambda$ on $\boO_\gamma(\C)$ by $\chi_\lambda(a)
=\sum a_n \lambda^n$.
\begin{prop}
\label{caracteres_Ogamma}
The characters of $\boO_\gamma(\C)$ are exactly the $\chi_\lambda$,
for $\lambda\in S^1$.
\end{prop}
\begin{proof}
This result is given by Rogozin in \cite{rogozin}, but there is a
(density) problem in his argument, for $\gamma \not \in \N$.
A correction is given in \cite{rogozin:corr}, and a more direct
argument can be found in \cite[Theorem 1.2.12]{frenk}.
\end{proof}

The following theorem has been thoroughly used in Section
\ref{section_resultat_abstrait}. It is a Wiener Lemma in the
algebra $\boO_\gamma(\C)$.
\begin{thm}
\label{wiener_S1}
Let $\boC$ be a Banach algebra, $\gamma>1$, and $A(z)=\sum_{n\in
\Z} A_n z^n \in \boO_\gamma(\boC)$. Assume that, for every $z\in
S^1$, $A(z)$ is an invertible element of $\boC$. Then $A$ is
invertible in the Banach algebra $\boO_\gamma(\boC)$.
\end{thm}
\begin{proof}
Gelfand's Theorem (\cite[Theorem 11.5 (c)]{rudin:funct_an2})
ensures that, if an element $a$ of a commutative Banach algebra
satisfies $\chi(a)\not=0$ for every character $\chi$, then $a$ is
invertible. With Proposition \ref{caracteres_Ogamma}, this gives
Theorem \ref{wiener_S1} for $\boO_\gamma(\C)$.

To handle the case of a general noncommutative
Banach algebra, we use Theorem 3 of
\cite{wiener_non_commutatif}.
\end{proof}

The same kind of Wiener Lemma holds in the algebra
$\boO_\gamma^+(\boC)$ (also defined in Section
\ref{subsection_wiener}):
\begin{thm}
\label{wiener_D}
Let $\boC$ be a Banach algebra, $\gamma>1$, and $A(z)=\sum_{n\in
\N} A_n z^n \in \boO_\gamma^+(\boC)$. Assume that, for every $z\in
\overline{\D}$, $A(z)$ is an invertible element of $\boC$. Then
$A$ is invertible in the Banach algebra $\boO_\gamma^+(\boC)$.
\end{thm}
\begin{proof}
This is the same proof as in Theorem \ref{wiener_S1} (but here,
the characters on $\boO_\gamma^+(\boC)$ are given by
$\chi_\lambda(a)=\sum_{n=0}^\infty \lambda^n a_n$ for $\lambda
\in \overline{\D}$).
\end{proof}

{\small
\bibliography{biblio}
\bibliographystyle{alpha}
}

\end{document}